\documentclass[a4paper,fleqn]{cas-sc}

\def\tsc#1{\csdef{#1}{\textsc{\lowercase{#1}}\xspace}}
\tsc{WGM}
\tsc{QE}
\tsc{EP}
\tsc{PMS}
\tsc{BEC}
\tsc{DE}

\usepackage[authoryear,longnamesfirst]{natbib}

\usepackage[T1]{fontenc}
\usepackage[latin9]{inputenc}
\usepackage{verbatim}
\usepackage{float}
\usepackage{amsmath}
\usepackage{amssymb}
\usepackage{stmaryrd}
\usepackage{graphicx}
\usepackage{booktabs}   

\makeatletter

\floatstyle{ruled}
\newfloat{algorithm}{tbp}{loa}
\providecommand{\algorithmname}{Algorithm}
\floatname{algorithm}{\protect\algorithmname}

\usepackage{babel}
\begin{document}

\let\WriteBookmarks\relax
\def\floatpagepagefraction{1}
\def\textpagefraction{.001}

\title{Efficient SN-like and PN-like Dynamic Low Rank methods for Thermal Radiative Transfer}
\shorttitle{SN-like and PN-like DLR methods for TRT}
\shortauthors{T. S. Haut et~al.}

\affiliation[1]{organization={Lawrence Livermore National Laboratory},
                addressline={7000 East Avenue}, 
                city={Livermore},
                postcode={94550}, 
                state={CA},
               country={United States}}

\affiliation[2]{organization={University of Innsbruck},
                addressline={ Innrain 52}, 
                city={ Innsbruck},
                postcode={6020}, 
               country={ Austria}}

\author[1]{T. Haut}
\author[1]{J. Loffeld}
\author[2]{L. Einkemmer}
\author[1]{P. Guthrey}
\author[2]{S. Brunner}
\author[1]{W. Schill}

\begin{abstract}
Dynamic Low Rank (DLR) methods are a promising way to reduce the computational
cost and memory footprint of the high-dimensional thermal radiative
transfer (TRT) equations. The TRT equations are a system of nonlinear
PDEs that model the energy exhchange between the material temperature
$T\left(\mathbf{x},t\right)$ and the radiation energy density $\psi\left(\mathbf{x},\boldsymbol{\Omega},t\right)$
at spatial point $\mathbf{x}\in\mathbb{R}^{3}$ and traveling in direction
$\boldsymbol{\Omega}\in\mathbb{S}^{2}$ at time $t$; due to their
high dimensionality, solving the TRT equations is often bottleneck
in multi-physics simulations. DLR methods represent the solution $\psi\left(\mathbf{x},\boldsymbol{\Omega},t\right)$
in terms of time-evolving SVD-like factors $W_{j}\left(\boldsymbol{\Omega},t\right)$
and $X_{j}\left(\mathbf{x},t\right)$ of angle and space. Although
previous work has explored DLR methods for TRT, most of the methods
have limitations that make them impractical for realistic scenarios  
and uncompetitive with current non-DLR production codes.

Here we develop new PN-like and SN-like Dynamic Low Rank (DLR) methods
for TRT. In the SN-like DLR method, we use the time-evolving angular
basis functions $W_{i}\left(\boldsymbol{\Omega},t\right)$ to select
time-evolving angles $\boldsymbol{\Omega}^{\left(i\right)}\left(t\right)\in\mathbb{S}^{2}$;
this DLR formulation enables us to use the highly optimized SN transport
sweep as our main computational kernel, and results in a practical
way of leveraging low-rank methods in production TRT codes. In contrast,
our PN-like DLR method uses an even-parity formulation and results
in positive-definite linear systems to solve for each time step.

We demonstrate the methods on several challenging, highly heterogenous problems
in two spatial dimensions $(4$D) that these DLR schemes can give
significant reduction in angular artifacts (``ray effects'') with
the same cost as gold-standard SN methods.
\end{abstract}

\maketitle

\section{Introduction}

Dynamic Low Rank (DLR) methods (\citep{koch2007dynamical}, \citep{EinkemmerLubich-2018},
\citep{Einkemmer-Review-2025}) are a promising approach for reducing
the computational cost and memory footprint of thermal radiative transfer.

We apply the DLR method to compute the time-dependent photon energy
distribution $\psi\left(\mathbf{x},\boldsymbol{\Omega},t\right)$
of the time-dependent thermal radiative transfer (TRT) equations
(see equations (\ref{eq:TRT. psi})-(\ref{eq:TRT, T})); here $\psi\left(\mathbf{x},\boldsymbol{\Omega},t\right)$
represents the radiation energy density at spatial point $\mathbf{x}\in\mathcal{X}\subset\mathbb{R}^{3}$,
traveling in direction $\boldsymbol{\Omega}\in\mathbb{S}^{2}$, at
time $t$. The DLR method represents the solution in terms of a dynamically
evolving low-rank representation \citep{EinkemmerLubich-2018} with
rank $r_{0}$,
\[
\psi\left(\mathbf{x},\boldsymbol{\Omega},t\right)=\sum_{i,j=1}^{r_{0}}X_{i}\left(\mathbf{x},t\right)S_{ij}\left(t\right)W_{j}\left(\boldsymbol{\Omega},t\right),
\]
where the $r_{0}$ spatial basis functions $X_{i}\left(\mathbf{x},t\right)$
and angular basis functions $W_{j}\left(\boldsymbol{\Omega},t\right)$
evolve in time to achieve a (nearly) optimally small equation residual
for a prescribed solution rank $r_{0}$. The memory footprint for
storing $X_{i}\left(\mathbf{x},t\right)$ and $W_{j}\left(\boldsymbol{\Omega},t\right)$
scales like $\mathcal{O}\left(r_{0}\left(N_{\mathbf{x}}+N_{\boldsymbol{\Omega}}\right)\right)$
and the computational cost scales like $\mathcal{O}\left(r_{0}^{2}\left(N_{\mathbf{x}}+N_{\boldsymbol{\Omega}}\right)\right)$,
with $N_{\mathbf{x}}$ and $N_{\boldsymbol{\Omega}}$ denoting the
number of spatial and angular unknowns needed to resolve the solution.
In contrast, direct numerical discretization schemes require $\mathcal{O}\left(N_{\mathbf{x}}N_{\boldsymbol{\Omega}}\right)$
memory and computational costs. Therefore, when $r_{0}^{2}\ll\max\left(N_{\mathbf{x}},N_{\boldsymbol{\Omega}}\right)$,
DLR offers significant computational savings.

There is a growing literature on applying DLR in the context of both
linear transport and (nonlinear) thermal radiative transfer. These
include conservative high-order low-order methods for radiative transfer
(cf. \citep{peng2021high}, \citep{Einkemmer2023c}, \citep{Einkemmer2021a},
\citep{Baumann2023}), as well as sweep-based methods on structured
meshes (cf. \citep{peng2020low}, \citep{peng2021low}, \citep{peng2023sweep}).
Many DLR schemes for kinetic equations such as TRT ensure that key
asymptotic-preserving properties hold (cf. \citep{Ceruti2023}, \citep{Frank2025}).

For TRT applications, there are still outstanding challenges for making
DLR practical in realistic scenarios compared to a gold-standard deterministic
method like the discrete-ordinates (SN) method. One major challenge
is that the number $N_{\boldsymbol{\Omega}}$ of angular collocation
points needed for acceptable physics fidelity in SN simulations can
be quite small in practice compared to the number of spatial unknowns
$N_{\mathbf{x}}$; in fact, radiation hydrodynamics simulations typically
use the hydrodynamics mesh to represent the spatial unknowns---which
can require  millions of spatial unknowns---but often only
require a few dozen to a few hundred angular directions. Second, there
are highly optimized, matrix-free transport solver methods---so called
transport sweeps (cf. \citep{mclendon2001finding,pautz2002algorithm,pautz2017parallel,plimpton2005parallel})---that
allow very efficient solutions of the SN equations (note that implicit
time-stepping schemes are needed for TRT, due to its stiffness). In
particular, SN methods do not require explicitly forming sparse matrices.
In contrast, a standard DLR scheme would involve forming and solving
the solution of an $r_{0}N_{\mathbf{x}}\times r_{0}N_{\mathbf{x}}$
non-symmetric sparse system. This makes it challenging to out-perform
SN methods in both memory and computational cost.

Inspired by the recent collocation-based DLR method in (cf. \citep{dektor2024interpolatory,dektor2025collocation},
\citep{Ghahremani2024}, \citep{Hossein2025}, and \citep{Zheng2025}),
we develop an SN-like DLR method. The scheme resembles a traditional
SN method, but with angular discretization points $\boldsymbol{\Omega}^{\left(i\right)}\left(t\right),i=1,\ldots,r_{0}$,
that are dynamically selected in time from the (evolving) DLR angular
basis functions $W_{i}\left(\boldsymbol{\Omega},t\right)$. This allows
us to use the highly optimized transport sweep solvers from a traditional
SN method as our primary computional kernel. In addition, we show
in both theory and practice that the traditional diffusion-synthetic
acceleration method (DSA) for preconditioning in the so-called diffusion
limit works without modification. As a result, we are able to use
the compression features of DLR in SN production codes with little
modification, and with essentially no additional computational overhead.
We note that a sweep-like method is developed in \citep{peng2023sweep},
but involves inverting $2^{\text{dim}}r_{0}\times2^{\text{dim}}r_{0}$
matrices for each spatial mesh element during the sweep and is limited
to structured meshes; in contrast, our approach involves inverting
only $2^{\text{dim}}\times2^{\text{dim}}$ matrices in each spatial
element for each of the $r_{0}$ directions, and can use state-of-the-art
transport sweeps for this purpose.

We also present an alternative, more traditional DLR method for TRT
that achieves efficiency in a different way.  As a classic DLR
method, it is based on Galerkin projection of the TRT equations on to angular
and spatial basis functions at the beginning of each time step.
This results in a PN-like DLR method, although in contrast to the
standard PN method---which employs spherical harmonics as angular
basis functions---the angular basis functions in the PN-like DLR
scheme evolve in time to (nearly) minimize the solution residual for
a fixed rank. Sweep solvers, which are SN-specific, are unavailable to PN methods,
and while there have been advances in applying multigrid solvers
to hyperbolic problems such as transport {\citep{Hanophy01112020, dargaville2025coarsening, dargaville2024air}, such approaches are not yet competitive
with sweep-based SN methods. However, in combination with DLR compression,
this PN-like DLR method attains efficiency through an even-parity reformulation
that transforms the first-order transport system into a pair
of second-order problems. This results in two positive-definite
$r_{0}N_{\mathbf{x}}\times r_{0}N_{\mathbf{x}}$ sparse matrix inversions
for each time step, for which efficient algebraic multigrid solvers can be used
for their solutions. In addition, the even-parity formulation
avoids the need for diffusion acceleration or iterations on a given time step,
as long as we ensure that the constant basis function in angle is contained
within our DLR angular basis.

Section~\ref{sec:Numerical-formulation-of} presents the TRT equations
and their reduction to a linear transport problem with effective scattering
via Newton linearization. Our PN-like and SN-like DLR
methods are presented in Sections~\ref{sec:PN like DLR scheme}-\ref{sec:SN like DLR scheme}.
Finally, Section~\ref{sec:Numerical-experiments} demonstrates these
DLR methods on several challenging test problems. The numerical experiments
show that the SN-like DLR method can give improved accuracy for the
same computational cost as a gold-standard SN method on a challenging
benchmark problem \citep{brunner2023family}.

\section{Numerical formulation of thermal radiative transfer (TRT) \label{sec:Numerical-formulation-of}}

We begin by describing the basic TRT system of equations.
The system describes the coupling between radiation energy density ($\psi\left(\mathbf{x},\boldsymbol{\Omega},t\right)$)
and material temperature ($T$).  $T$ is lower dimensional and is not factorized by DLR,
so we also discuss the linearized time stepping approach for the temperature equation
that is used throughout by both standard SN and DLR methods. The DLR factorization
applies to $\psi\left(\mathbf{x},\boldsymbol{\Omega},t\right)$, so we leave the description
of discretization of the radiation energy density equation within each DLR method to their respective sections.

\subsection{TRT equations}

The frequency-averaged (i.e., grey) TRT equations for the material
temperature $T\left(\mathbf{x},t\right)$ and the photon energy distribution
function $\psi\left(\mathbf{x},\boldsymbol{\Omega},t\right)$ are
given by \citep{castor2004radiation}
\begin{equation}
\frac{1}{c}\frac{\partial\psi}{\partial t}+\boldsymbol{\Omega}\cdot\nabla_{\mathbf{x}}\psi+\sigma\left(T\right)\psi=\sigma\left(T\right)B\left(T\right),\,\,\,\,\mathbf{x}\in\mathcal{X},\label{eq:TRT. psi}
\end{equation}
\begin{equation}
\psi\left(\mathbf{x},\boldsymbol{\Omega},t\right)=0,\,\,\,\,\boldsymbol{\Omega}\cdot\mathbf{n}\left(\mathbf{x}\right)\leq0,\,\,\,\,\mathbf{x}\in\mathcal{\partial X},\label{eq:TRT, inflow}
\end{equation}
\begin{equation}
\rho c_{v}\left(T\right)\frac{\partial T}{\partial t}=\sigma\left(T\right)\left(\varphi-4\pi B\left(T\right)\right).\label{eq:TRT, T}
\end{equation}
Here $\rho$ denotes the material density, $c_{v}\left(T\right)$
denotes the temperature-dependent specific heat, $\sigma\left(T\right)$
denotes the temperature-dependent absorption opacity, $\varphi$ denotes
the angularly averaged intensity,
\[
\varphi\left(\mathbf{x}\right)=\int_{\mathbb{S}^{2}}\psi\left(\mathbf{x},\boldsymbol{\Omega}\right)d\boldsymbol{\Omega}.
\]
For simplicity, we only consider vacuum boundary conditions (\ref{eq:TRT, inflow})
in this paper; however, the presented DLR methods extend trivially
to more general boundary conditions.

\subsection{Backward Euler discretization of the temperature equation and its
linearization}

We discretize equation (\ref{eq:TRT, T}) as follows

\begin{equation}
\rho c_{v}\left(T_{0}\right)\left(T-T_{0}\right)=\Delta t\sigma\left(T_{0}\right)\left(\varphi-4\pi B\left(T\right)\right).\label{eq:T, backward Euler}
\end{equation}
As is standard, we treat the temperature dependence of the opacity
and specific heat explicitly, and treat the stiff emission term $B\left(T\right)$
implicitly.

As is standard, we obtain an approximation of the new temperature
$T$ in equation (\ref{eq:T, backward Euler}) via a Newton linearization
\citep{aarseth1985multiple} about the current temperature $T_{0}$.
This linearization yields a local $\mathcal{O}\left(\Delta t^{2}\right)$
time-stepping error, and results in an overall time-stepping scheme---for
both DLR and standard SN methods---that in practice is as stable
as the fully implicit backward Euler scheme.

In particular, a standard calculation yields the linear transport
equation with pseudo-scattering over the time step $t\in\left[t_{0},t_{0}+\Delta t\right]$,
\begin{equation}
\frac{1}{c}\frac{\partial\psi}{\partial t}+\boldsymbol{\Omega}\cdot\nabla_{\mathbf{x}}\psi+\sigma\psi=\frac{1}{4\pi}\sigma_{s}\left(T_{0}\right)\varphi+q\left(T_{0}\right),\label{eq:linear transport with pseudoscattering}
\end{equation}
where the pseudo-scattering opacity $\sigma_{s}\left(T_{0}\right)$
and the source term $q\left(T_{0}\right)$ are defined as
\begin{equation}
\sigma_{s}\left(T_{0}\right)=\left(\frac{4\pi\Delta t\sigma\left(T_{0}\right)\frac{\partial B}{\partial T}\left(T_{0}\right)}{\rho c_{v}\left(T_{0}\right)+4\pi\Delta t\sigma\left(T_{0}\right)\frac{\partial B}{\partial T}\left(T_{0}\right)}\right)\sigma\left(T_{0}\right),\label{eq:effective scattering}
\end{equation}
\begin{equation}
q\left(T_{0}\right)=\sigma\left(T_{0}\right)B\left(T_{0}\right)-\sigma\frac{\partial B}{\partial T}\left(T_{0}\right)\frac{4\pi\Delta t\sigma\left(T_{0}\right)B\left(T_{0}\right)+\rho c_{v}\left(T_{0}\right)T_{0}}{\rho c_{v}+4\pi\sigma\left(T_{0}\right)\Delta t\frac{\partial B}{\partial T}\left(T_{0}\right)}.\label{eq:effective source from Newton}
\end{equation}
Once we solve for $\psi$, we can update our approximation to temperature
at the next time step,
\begin{equation}
T\approx T_{0}+\frac{\Delta t\sigma\left(T_{0}\right)\left(\varphi-4\pi B\left(T_{0}\right)\right)-\rho c_{v}\left(T_{0}\right)T_{0}}{\rho c_{v}\left(T_{0}\right)+4\pi\Delta t\sigma\left(T_{0}\right)\frac{\partial B}{\partial T}\left(T_{0}\right)}.\label{eq:temperature update}
\end{equation}

\section{A PN-like DLR scheme using an even-parity formulation of transport\label{sec:PN like DLR scheme}}

Here we develop a classical DLR based TRT scheme, with the twist that we employ an even-parity transformation
that allows the system to be solved efficiently using algebraic multigrid solvers.  The result is a PN-like DLR formulation
that is efficiently solvable and exhibits no ray effects

We approximate the solution using a time-evolving SVD-like low-rank
representation

\begin{equation}
\psi\left(\mathbf{x},\boldsymbol{\Omega},t\right)=\sum_{i,j=1}^{r_{0}}X_{i}\left(\mathbf{x},t\right)S_{ij}\left(t\right)W_{j}\left(\boldsymbol{\Omega},t\right)=\boldsymbol{W}^{\text{T}}\left(\boldsymbol{\Omega},t\right)S\left(t\right)^{\text{T}}\mathbf{X}\left(\mathbf{x},t\right),\label{eq:low rank}
\end{equation}
where
\[
\mathbf{X}\left(\mathbf{x},t\right)=\left(\begin{array}{c}
X_{1}\left(\mathbf{x},t\right)\\
\vdots\\
X_{r_{0}}\left(\mathbf{x},t\right)
\end{array}\right),\,\,\,\,\,\mathbf{W}\left(\boldsymbol{\Omega},t\right)=\left(\begin{array}{c}
W_{1}\left(\boldsymbol{\Omega},t\right)\\
\vdots\\
W_{r_{0}}\left(\boldsymbol{\Omega},t\right)
\end{array}\right).
\]
Here the functions $X_{i}$, $i=1,\ldots,r_0$, are orthonormal,
\begin{equation}
\left\langle X_{i}\left(\cdot,t\right),X_{j}\left(\cdot,t\right)\right\rangle _{\mathcal{X}}=\int_{\mathcal{X}}X_{i}\left(\mathbf{x},t\right)X_{j}\left(\mathbf{x},t\right)d\mathbf{x}=\delta_{i,j}.\label{eq:orthonormality of X}
\end{equation}
Similarly, the functions $W_{j}$ are orthonormal,
\begin{equation}
\left\langle W_{i}\left(\cdot,t\right),W_{j}\left(\cdot,t\right)\right\rangle _{\mathcal{X}}=\int_{\mathbb{S}^{2}}W_{i}\left(\boldsymbol{\Omega},t\right)W_{j}\left(\boldsymbol{\Omega},t\right)d\boldsymbol{\Omega}=\delta_{i,j}.\label{eq::orthonormality of W}
\end{equation}
We frequently express the low-rank representation (\ref{eq:low rank})
in the equivalent forms
\begin{align}
\psi\left(\mathbf{x},\boldsymbol{\Omega},t\right) & =\sum_{i,j=1}^{r_{0}}X_{i}\left(\mathbf{x},t\right)S_{ij}\left(t\right)W_{j}\left(\boldsymbol{\Omega},t\right)\nonumber \\
 & =\sum_{j=1}^{r_{0}}K_{j}\left(\mathbf{x},t\right)W_{j}\left(\boldsymbol{\Omega},t\right)\label{eq:DLR decomposition, K and W}\\
 & =\sum_{i=1}^{r_{0}}X_{i}\left(\mathbf{x},t\right)L_{i}\left(\boldsymbol{\Omega},t\right),\label{eq:DLR decomposition, X and L}
\end{align}
where
\begin{equation}
K_{j}\left(\mathbf{x},t\right)=\sum_{i=1}^{r_{0}}X_{i}\left(\mathbf{x},t\right)S_{ij}\left(t\right),\,\,\,\,\,\,\,\,\,\,L_{i}\left(\boldsymbol{\Omega},t\right)=\sum_{j=1}^{r_{0}}S_{ij}\left(t\right)W_{j}\left(\boldsymbol{\Omega},t\right).\label{eq:K and L}
\end{equation}

Write the transport equation (\ref{eq:linear transport with pseudoscattering})
in the abstract form
\begin{equation}
\frac{\partial\psi}{\partial t}=\mathcal{L}\psi+q.\label{eq:transport, abstract}
\end{equation}

Differentiating equations (\ref{eq:DLR decomposition, K and W}) and
(\ref{eq:DLR decomposition, X and L}) with respect to time $t$ and
using the orthogonality constraints (\ref{eq:orthonormality of X})-(\ref{eq::orthonormality of W}),
it is straightforward to show (see Appendix~\ref{subsec:Appendix A})
that 
\begin{equation}
\dot{L}_{j}=\left\langle \frac{\partial\psi}{\partial t},X_{j}\right\rangle _{\mathcal{X}},\,\,\,\,\,\dot{K}_{j}=\left\langle \frac{\partial\psi}{\partial t},W_{j}\right\rangle _{\mathbb{S}^{2}},\,\,\,\,\,\dot{S}_{ij}=\left\langle \frac{\partial\psi}{\partial t},X_{i}W_{j}\right\rangle _{\mathcal{X}\times\mathbb{S}^{2}},\label{eq:L dot, K dot, S dot}
\end{equation}
and 
\begin{equation}
\frac{\partial\psi}{\partial t}=\sum_{j=1}^{r_0}\dot{K}_{j}W_{j}+\sum_{i=1i}^{r_0}\dot{L}_{i}X_{i}-\sum_{i,j=1}^{r_0}X_{i}\dot{S}_{ij}W_{j}.\label{eq:dpsi_dt}
\end{equation}
Combining equations (\ref{eq:transport, abstract}) and (\ref{eq:L dot, K dot, S dot}),
\begin{equation}
\dot{L}_{j}=\sum_{i=1i}^{r_0}\left\langle \mathcal{L}\left(X_{i}L_{i}\right)+q,X_{j}\right\rangle _{\mathcal{X}},\label{eq:L_dot}
\end{equation}
\begin{equation}
\dot{K}_{j}=\sum_{i=1}^{r_0}\left\langle \mathcal{L}\left(K_{i}W_{i}\right)+q,W_{j}\right\rangle _{\mathbb{S}^{2}}.\label{eq:K_dot}
\end{equation}
\begin{equation}
\dot{S}_{ij}=\sum_{i',j'=1}^{r_0}S_{i'j'}\left\langle \mathcal{L}\left(X_{i'}W_{j'}\right)+q,X_{i}W_{j}\right\rangle _{\mathcal{X}\times\mathbb{S}^{2}},\label{eq:S_dot}
\end{equation}

We discretize equations (\ref{eq:L_dot})-(\ref{eq:S_dot}) in time
via backward Euler and lag the DLR basis functions in time $X_{j}$
and $W_{j}$. In particular, define $X_{0,i}\left(\mathbf{x}\right)=X_{i}\left(\mathbf{x},t_{0}\right)$
and $W_{0,j}\left(\boldsymbol{\Omega}\right)=W_{j}\left(\boldsymbol{\Omega},t_{0}\right)$.
Then
\begin{equation}
\frac{L_{j}-L_{0,j}}{\Delta t}=\sum_{i=1}^{r_{0}}\left\langle \mathcal{L}\left(X_{0,i}L_{i}\right)+q,X_{0,j}\right\rangle _{\mathcal{X}},\label{eq:DLR Galerkin, L step}
\end{equation}
\begin{equation}
\frac{K_{j}-K_{0,j}}{\Delta t}=\sum_{i=1}^{r_{0}}\left\langle \mathcal{L}\left(W_{0,i}K_{i}\right)+q,W_{0,j}\right\rangle _{\mathbb{S}^{2}},\label{eq:DLR Galerkin, K step}
\end{equation}
\begin{equation}
\frac{S_{ij}-S_{0,ij}}{\Delta t}=\sum_{i',j'=1}^{r_{0}}S_{i'j'}\left\langle \mathcal{L}\left(X_{0,i'}W_{0,j'}\right)+q,X_{0,i}W_{0,j}\right\rangle _{\mathcal{X}\times\mathbb{S}^{2}}.\label{eq:DLR Galerkin, S step}
\end{equation}
We supplement the above equations with equation (\ref{eq:T, backward Euler})
for the updated temperature $T$.

Once we solve the equations (\ref{eq:T, backward Euler}) and (\ref{eq:DLR Galerkin, L step})-(\ref{eq:DLR Galerkin, S step})
for $T$, $K_{j}$, $L_{j}$, and $S_{ij}$ at the next time $t=t_{0}+\Delta t$,
we then compute the updated solution via equation (\ref{eq:dpsi_dt}),
\begin{align}
\psi & =\psi_{0}+\Delta t\frac{\partial\psi}{\partial t}\nonumber \\
 & =\psi_{0}+\sum_{j=1}^{r_{0}}\left(K_{j}-K_{0,j}\right)W_{j}+\sum_{i=1}^{r_{0}}\left(L_{i}-L_{0,i}\right)X_{i}-\sum_{i,j=1}^{r_{0}}X_{i}\left(S_{ij}-S_{0,ij}\right)W_{j}.\label{eq:psi updated, rank 2r}
\end{align}

The representation (\ref{eq:psi updated, rank 2r}) for $\psi$ at
the next time step increases the rank from $r_{0}$ to $2r_{0}$;
we perform a truncated SVD to reduce the representation back down
to rank $r_{0}$. In more detail, 
\[
\psi\left(\mathbf{x},\boldsymbol{\Omega}\right)=\mathbf{Y}\left(\mathbf{x}\right)^{\text{T}}\mathbf{Z}\left(\boldsymbol{\Omega}\right),
\]
where, e.g.,
\[
Y_{j}=\begin{cases}
\left(K_{j}-K_{0,j}\right)-\sum_{i=1}^{r_{0}}X_{i}\left(S_{ij}-S_{0,ij}\right), & 1\leq j\leq r_{0},\\
X_{i}, & r_{0}+1\leq j\leq2r_{0}.
\end{cases}
\]
We can orthogonalize the basis functions $Y_{j}\left(\mathbf{x}\right)$
via stabilized Gram Schmidt, $\mathbf{Y}\left(\mathbf{x}\right)=R_{Y}\mathbf{U}\left(\mathbf{x}\right)$.
Similarly, $\mathbf{Z}\left(\boldsymbol{\Omega}\right)=R_{Z}\mathbf{V}\left(\boldsymbol{\Omega}\right)$.
Then
\[
\psi\left(\mathbf{x},\boldsymbol{\Omega}\right)=\mathbf{Y}\left(\mathbf{x}\right)^{\text{T}}\mathbf{Z}\left(\boldsymbol{\Omega}\right)=\mathbf{U}\left(\mathbf{x}\right)^{\text{T}}\left(R_{Y}^{\text{T}}R_{Z}\right)\mathbf{V}\left(\boldsymbol{\Omega}\right).
\]
We then employ a truncated SVD on the $2r_{0}\times2r_{0}$ matrix
$R_{Y}^{\text{T}}R_{Z}$ to round back down to a rank $r$ approximation
for $\psi\left(\mathbf{x},\boldsymbol{\Omega}\right)$.

We now discuss the solution of equations (\ref{eq:DLR Galerkin, L step})-(\ref{eq:DLR Galerkin, S step})
in more detail in Sections~\ref{subsec:Even-parity formulation}-\ref{subsec:S step for PN-lie DLR method}.

\subsection{Even parity formulation of transport for the K step\label{subsec:Even-parity formulation}}

Fast transport sweeps are unavailable to PN formulations, and multigrid solvers do not yet give competitive
efficiency when inverting the matrices that come from direct discretization of equation (\ref{eq:linear transport with pseudoscattering}).
To deal with this, we use the even parity approach from transport to reformulate the first-order system into a pair of second-order systems
(cf. \cite{lewis1977progress}) that are amenable to fast algebraic multigrid solvers.  This trades one slow solve for two fast ones,
and in combination with DLR compression gives an efficient solver.  Here we briefly describe the even parity formulation for standard
transport and then adapt its use to DLR in the next section.

Consider the backward Euler discretization of the steady-state equation
(\ref{eq:linear transport with pseudoscattering}),
\begin{equation}
\boldsymbol{\Omega}\cdot\nabla_{\mathbf{x}}\psi+\sigma_{t}\psi-\frac{1}{4\pi}\sigma_{s}\varphi=Q.\label{eq:transport}
\end{equation}
First we decompose the solution into even and odd components in angle,
\[
\psi_{\pm}\left(\mathbf{x},\boldsymbol{\Omega}\right)=\frac{1}{2}\left(\psi\left(\mathbf{x},\boldsymbol{\Omega}\right)\pm\psi\left(\mathbf{x},-\boldsymbol{\Omega}\right)\right).
\]
Then it can be shown (cf. \citep{lewis1977progress}) that $\psi_{-}$
satisfies 
\begin{equation}
-\boldsymbol{\Omega}\cdot\nabla_{\mathbf{x}}\left(\frac{1}{\sigma_{t}}\boldsymbol{\Omega}\cdot\nabla\psi_{-}\right)+\sigma_{t}\psi_{-}=Q_{-}-\boldsymbol{\Omega}\cdot\nabla_{\mathbf{x}}\left(\frac{1}{\sigma_{t}}Q_{+}\right),\,\,\,\,\mathbf{x}\in\mathcal{X},\label{eq:even parity, --1}
\end{equation}
\begin{equation}
\left|\mathbf{n}\cdot\boldsymbol{\Omega}\right|\psi_{-}+\mathbf{n}\cdot\boldsymbol{\Omega}\left(\frac{1}{\sigma_{t}}\boldsymbol{\Omega}\cdot\nabla_{\mathbf{x}}\psi_{-}\right)=\mathbf{n}\cdot\boldsymbol{\Omega}\frac{1}{\sigma_{t}}Q_{+},\,\,\,\,\,\mathbf{x}\in\partial\mathcal{X}.\label{eq:even parity -, boundary}
\end{equation}
Similarly,
\begin{equation}
-\boldsymbol{\Omega}\cdot\nabla_{\mathbf{x}}\left(\frac{1}{\sigma_{t}}\boldsymbol{\Omega}\cdot\nabla_{\mathbf{x}}\psi_{+}\right)+\sigma_{t}\psi_{+}-\frac{\sigma_{s}}{4\pi}\varphi_{+}=Q_{+}-\boldsymbol{\Omega}\cdot\nabla_{\mathbf{x}}\left(\frac{1}{\sigma_{t}}Q_{-}\right),\,\,\,\,\mathbf{x}\in\mathcal{X},\label{eq:even parity, +}
\end{equation}
\begin{equation}
\left|\mathbf{n}\cdot\boldsymbol{\Omega}\right|\psi_{+}+\mathbf{n}\cdot\boldsymbol{\Omega}\left(\frac{1}{\sigma_{t}}\boldsymbol{\Omega}\cdot\nabla_{\mathbf{x}}\psi_{+}\right)=\mathbf{n}\cdot\boldsymbol{\Omega}\frac{1}{\sigma_{t}}Q_{-},\,\,\,\,\,\mathbf{x}\in\partial\mathcal{X}.\label{eq:even parity boundary +, boundary}
\end{equation}

\subsection{K step for PN-like DLR method\label{subsec:K step for PN-lie DLR method}}

We now explicitly compute the matrix system for equations (\ref{eq:DLR Galerkin, K step}),
starting from the even-parity formulation detailed in Section~\ref{subsec:Even-parity formulation}.

Importantly, we always enforce that $W_{0,0}\left(\boldsymbol{\Omega}\right)=1/\sqrt{4\pi}$.
This ensures that Galerkin projection for the K step can be solved
efficiently without iteration. We also use even and odd constructions
of the angular basis functions,
\begin{equation}
W_{0,j}^{+}\left(\boldsymbol{\Omega}\right)=\frac{W_{0,j}\left(\boldsymbol{\Omega}\right)+W_{0,j}\left(-\boldsymbol{\Omega}\right)}{2},\,\,\,\,\,\,\,\,W_{0,j}^{-}\left(\boldsymbol{\Omega}\right)=\frac{W_{0,j}\left(\boldsymbol{\Omega}\right)-W_{0,j}\left(-\boldsymbol{\Omega}\right)}{2},\label{eq:W even and odd}
\end{equation}
when we construct the K step equations.

To do so, suppose that
\[
\psi_{\pm}\left(\mathbf{x},\boldsymbol{\Omega}\right)=\sum_{j=1}^{r_{0}}W_{0,j}^{\pm}\left(\boldsymbol{\Omega}\right)K_{j}^{\pm}\left(\mathbf{x}\right).
\]
Define 
\[
A_{j,j'}^{\pm}=\left\langle \boldsymbol{\Omega}\boldsymbol{\Omega}^{\text{T}}W_{0,j}^{\pm},W_{0,j'}^{\pm}\right\rangle _{\mathbb{S}^{2}},\,\,\,\,\,\,\boldsymbol{\kappa}_{j,j'}=\left\langle \boldsymbol{\Omega}W_{0,j}^{\pm},W_{0,j'}^{\pm}\right\rangle _{\mathbb{S}^{2}},\,\,\,\,\,\,\beta_{j,j'}\left(\mathbf{x}\right)=\left\langle \left|\mathbf{n}\left(\mathbf{x}\right)\cdot\boldsymbol{\Omega}\right|W_{0,j}^{\pm},W_{0,j'}^{\pm}\right\rangle _{\mathbb{S}^{2}}.
\]
Then projecting equations (\ref{eq:even parity, +}) and (\ref{eq:even parity boundary +, boundary})
on to the angular basis functions $W_{0,j}^{\pm}\left(\boldsymbol{\Omega}\right)$,
a direct calculation shows that, $j=1,\ldots,r_{0}$,
\begin{equation}
-\sum_{j=1}^{r_{0}}A_{j,j'}^{+}\nabla_{\mathbf{x}}\left(\frac{1}{\sigma_{t}}\nabla_{\mathbf{x}} K_{j}^{+}\right)+\sigma_{t}K_{j'}^{+}+\delta_{j',0}\sqrt{\frac{1}{4\pi}}\sigma_{s}K_{0}^{+}=Q_{j'}^{+}-\sum_{j=1}^{r_{0}}\nabla_{\mathbf{x}}\cdot\left(\frac{1}{\sigma_{t}}Q_{j'}^{-}\boldsymbol{\kappa}_{j,j'}^{+}\right),\,\,\,\,\mathbf{x}\in\mathcal{X},\label{eq:K step DLR Galerkin}
\end{equation}
with boundary conditions 
\begin{equation}
\sum_{j=1}^{r_{0}}\beta_{j,j'}^{+}K_{j}^{+}+\sum_{j=1}^{r_{0}}\mathbf{n}\cdot\frac{1}{\sigma_{t}}\left(A_{j,j'}^{+}\nabla_{\mathbf{x}}K_{j}^{+}\right)=\sum_{j=1}^{r_{0}}\mathbf{n}\cdot\frac{1}{\sigma_{t}}\left(\boldsymbol{\kappa}_{j,j'}^{+}Q_{j}^{-}\right),\,\,\,\,\,\mathbf{x}\in\partial\mathcal{X}.\label{eq:K step DLR Galerkin, boundary condition}
\end{equation}
An analogous system holds for $K_{j}^{-}$, $j=1,\ldots,r$:
\begin{equation}
-\sum_{j=1}^{r_{0}}A_{j,j'}^{-}\nabla_{\mathbf{x}}\left(\frac{1}{\sigma_{t}}\nabla_{\mathbf{x}} K_{j}^{-}\right)+\sigma_{t}K_{j'}^{-}=Q_{j'}^{-}-\sum_{j=1}^{r_{0}}\nabla_{\mathbf{x}}\cdot\left(\frac{1}{\sigma_{t}}Q_{j'}^{+}\boldsymbol{\kappa}_{j,j'}^{-}\right),\,\,\,\,\mathbf{x}\in\mathcal{X},\label{eq:plus, K step DLR Galerkin}
\end{equation}
with boundary conditions 
\begin{equation}
\sum_{j=1}^{r_{0}}\beta_{j,j'}^{-}K_{j}^{-}+\sum_{j=1}^{r_{0}}\mathbf{n}\cdot\frac{1}{\sigma_{t}}\left(A_{j,j'}^{-}\nabla_{\mathbf{x}}K_{j}^{-}\right)=\sum_{j=1}^{r_{0}}\mathbf{n}\cdot\frac{1}{\sigma_{t}}\left(\boldsymbol{\kappa}_{j,j'}^{-}Q_{j}^{+}\right),\,\,\,\,\,\mathbf{x}\in\partial\mathcal{X}.\label{eq:plus, K step DLR Galerkin, boundary condition}
\end{equation}

Note that the equations (\ref{eq:plus, K step DLR Galerkin})-(\ref{eq:plus, K step DLR Galerkin, boundary condition})
are positive-definite and their space-discretized version are efficiently
solved with algebraic multigrid methods.

\subsection{L step for PN-like DLR method\label{subsec:L step for PN-lie DLR method}}

Define
\[
B_{i,j}\left(\boldsymbol{\Omega}\right)=\left\langle \nabla_{\mathbf{x}}X_{0,j},X_{0,i}\right\rangle _{\mathcal{X}}\cdot\boldsymbol{\Omega}+\sum_{j}\left\langle \sigma_{t}X_{0,j},X_{0,i}\right\rangle _{\mathcal{X}},
\]
\[
\left(M_{s}\right)_{i,j}=\left\langle \sigma_{s}X_{0,j},X_{0,i}\right\rangle _{\mathcal{X}}.
\]
Then projecting equation (\ref{eq:transport}) on to the $r_{0}$
spatial basis functions $X_{0,j}\left(\mathbf{x}\right)$, a straight-forward
calculation shows that
\begin{equation}
\sum_{j=1}^{r_0}B_{i,j}\left(\boldsymbol{\Omega}\right)L_{j}\left(\boldsymbol{\Omega}\right)-\frac{1}{4\pi}\sum_{j=1}^{r_0}\left(M_{s}\right)_{i,j}\int_{\mathbb{S}^{2}}L_{j}d\boldsymbol{\Omega}'=q_{0,i}\left(\boldsymbol{\Omega}\right).\label{eq:L system}
\end{equation}

To efficiently solve the L system (\ref{eq:L system}), define

\[
\mathbf{L}\left(\boldsymbol{\Omega}\right)=\left(\begin{array}{c}
L_{1}\left(\boldsymbol{\Omega}\right)\\
\vdots\\
L_{r}\left(\boldsymbol{\Omega}\right)
\end{array}\right),\,\,\,\,\,\mathbf{q}\left(\boldsymbol{\Omega}\right)=\left(\begin{array}{c}
q_{0.1}\left(\boldsymbol{\Omega}\right)\\
\vdots\\
q_{0,r}\left(\boldsymbol{\Omega}\right)
\end{array}\right),\,\,\,\,\,\overline{\mathbf{L}}=\left(\begin{array}{c}
\int_{\mathbb{S}^{2}}L_{1}\left(\boldsymbol{\Omega}'\right)d\boldsymbol{\Omega}'\\
\vdots\\
\int_{\mathbb{S}^{2}}L_{r}\left(\boldsymbol{\Omega}'\right)d\boldsymbol{\Omega}'
\end{array}\right).
\]
Then 
\begin{equation}
B\left(\boldsymbol{\Omega}\right)L\left(\boldsymbol{\Omega}\right)-\frac{1}{4\pi}M_{s}\overline{\mathbf{L}}=\mathbf{q}\left(\boldsymbol{\Omega}\right).\label{eq:L matrix system}
\end{equation}
Define
\[
\overline{C}=\int_{\mathbb{S}^{2}}\left(B\left(\boldsymbol{\Omega}\right)\right)^{-1}d\Omega,\,\,\,\,\,\overline{\boldsymbol{\zeta}}=\int_{\mathbb{S}^{2}}B\left(\boldsymbol{\Omega}'\right)^{-1}\mathbf{q}\left(\boldsymbol{\Omega}'\right)d\Omega'.
\]
Then
\[
\overline{\mathbf{L}}-\overline{C}\frac{1}{4\pi}M_{s}\overline{\mathbf{L}}=\overline{\boldsymbol{\zeta}}.
\]
Once we solve for $\overline{\mathbf{L}}$, we solve equation (\ref{eq:L matrix system})
for each angle $\boldsymbol{\Omega}$ independently.

\subsection{S step for PN-like DLR method\label{subsec:S step for PN-lie DLR method}}

Define 
\[
\mathbf{b}_{i,i'}=\left\langle \nabla_{\mathbf{x}}X_{0,i},X_{0,i'}\right\rangle _{\mathcal{X}},\,\,\,\,\,\boldsymbol{\kappa}_{j,j'}=\left\langle \Omega W_{0,j},W_{0,j'}\right\rangle _{\mathbb{S}^{2}}.
\]
Then projecting equation (\ref{eq:transport}) on to the $r_{0}^{2}$
basis functions $X_{0,i}\left(\mathbf{x}\right)W_{0,j}\left(\boldsymbol{\Omega}\right)$,

\[
\sum_{i,j=1}^{r_{0}}S_{ij}\mathbf{b}_{i,i'}\boldsymbol{\kappa}_{j,j'}+\sum_{i,j}\left(M_{t}\right)_{i,i'}S_{ij}\delta_{j,j'}-\sum_{i,j=1}^{r_{0}}\frac{1}{4\pi}S_{ij}\left(M_{s}\right)_{i,i'}\overline{W}_{j'}\overline{W}_{j}=q_{i',j'}.
\]
Here
\[
q_{i',j'}=S_{0,i',j'}+\left\langle q,X_{0,i'}W_{0,j'}\right\rangle _{\mathcal{X}\times\mathbb{S}^{2}}.
\]

\section{An SN-like DLR scheme\label{sec:SN like DLR scheme}}

We present an SN-like DLR method for solving the TRT equations. The
key advantage of this scheme is that we can leverage highly efficient
sweep-based methods

We note that the sweep-based method in \citep{peng2023sweep} requires
the inversion of $r_{0}2^{\text{dim}}\times r_{0}2^{\text{dim}}$
matrices in each mesh element; our collocation-based scheme instead
requires the inversion of matrices of size $2^{\text{dim}}\times2^{\text{dim}}$
in each mesh element for rank $r_{0}$ worth of angles. In practice,
this makes this version of the sweeping scheme competetive with SN
methods even for small SN orders. In addition, our collocation-based
scheme mitigates, but does not eliminate ray effects, as demonstrated
by the hohlraum numerical experiments in Section~\ref{subsec:Hohlraum problem};
in such cases, there can be major advantages to the scheme in \citep{peng2023sweep},
as well as the PN-like DLR method presented in this section.

\subsection{Computing the scalar flux from collocation angles \label{subsec:Computing-the-scalar}}

Suppose that we have approximations to $\psi\left(\mathbf{x},\boldsymbol{\Omega}_{i}\left(t\right),t\right)$
at specific angles $\boldsymbol{\Omega}^{\left(i\right)}\left(t\right)$,
$i=1,\ldots,r_{0}$. We want to approximate $\psi\left(\mathbf{x},\boldsymbol{\Omega},t\right)$
in terms of the orthonormal angular functions $W_{j}\left(\boldsymbol{\Omega},t\right)$,
\begin{equation}
\psi\left(\mathbf{x},\boldsymbol{\Omega},t\right)\approx\sum_{j=1}^{r_{0}}K_{j}\left(\mathbf{x},t\right)W_{j}\left(\boldsymbol{\Omega},t\right).\label{eq:K_j W_j}
\end{equation}
To do so, we solve the linear system
\begin{align*}
\psi\left(\mathbf{x},\boldsymbol{\Omega}^{\left(i\right)}\left(t\right),t\right) & =\sum_{j=1}^{r_{0}}K_{j}\left(\mathbf{x},t\right)W_{j}\left(\boldsymbol{\Omega}^{\left(i\right)}\left(t\right),t\right),\,\,\,\,\,\,\,i=1,\ldots,r_{0}.
\end{align*}
Define $\hat{W}_{ij}\left(t\right)=W_{j}\left(\boldsymbol{\Omega}^{\left(i\right)}\left(t\right),t\right)$,
and
\[
\mathbf{K}\left(\mathbf{x},t\right)=\left(\begin{array}{c}
K_{1}\left(\mathbf{x},t\right)\\
\vdots\\
K_{r_{0}}\left(\mathbf{x},t\right)
\end{array}\right),\,\,\,\,\,\,\boldsymbol{\Psi}\left(\mathbf{x},t\right)=\left(\begin{array}{c}
\psi\left(\mathbf{x},\boldsymbol{\Omega}^{\left(1\right)}\left(t\right),t\right)\\
\vdots\\
\psi\left(\mathbf{x},\boldsymbol{\Omega}^{\left(r_{0}\right)}\left(t\right),t\right)
\end{array}\right).
\]
Then
\[
\boldsymbol{\Psi}\left(\mathbf{x},t\right)=\hat{W}\left(t\right)\mathbf{K}\left(\mathbf{x},t\right).
\]
It follows that we can reconstruct 
\[
\mathbf{K}\left(\mathbf{x},t\right)=\hat{W}\left(t\right)^{-1}\boldsymbol{\Psi}\left(\mathbf{x},t\right).
\]
From equation (\ref{eq:K_j W_j}), we also have that
\begin{align*}
\int_{\mathbb{S}^{2}}\psi\left(\mathbf{x},\boldsymbol{\Omega},t\right)d\boldsymbol{\Omega} & \approx\sum_{j=1}^{r_{0}}K_{j}\left(\mathbf{x},t\right)\int_{\mathbb{S}^{2}}W_{j}\left(\boldsymbol{\Omega},t\right)d\boldsymbol{\Omega}\\
 & =\boldsymbol{\beta}^{\text{T}}\left(\hat{W}^{-1}\boldsymbol{\Psi}\left(\mathbf{x},t\right)\right),
\end{align*}
where
\[
\beta_{j}=\int_{\mathbb{S}^{2}}W_{j}\left(\boldsymbol{\Omega}\right)d\boldsymbol{\Omega},\,\,\,\,\,\,j=1,\ldots,r_{0}.
\]

\subsection{Discrete Empirical Interpolation Method (DEIM) \label{sec:DEIM}}

The DEIM algorithm proceeds iteratively in a manner entirely analogous
to Gram-Schmidt orthogonalization.

On the first step, we choose the first interpolation point $\boldsymbol{\Omega}_{1}\in\mathbb{S}^{2}$
that maximizes $\left|W_{1}\left(\boldsymbol{\Omega}\right)\right|$,
\[
\boldsymbol{\Omega}_{1}=\text{argmax}_{\boldsymbol{\Omega}}\left|W_{1}\left(\boldsymbol{\Omega}\right)\right|.
\]

Now, suppose at the $\left(j-1\right)$st step we have constructed
interpolation points $\boldsymbol{\Omega}_{1},\ldots,\boldsymbol{\Omega}_{j-1}$.
Define the interpolatory projection
\[
\left(P^{\left(j-1\right)}W_{j}\right)\left(\boldsymbol{\Omega}\right)=\sum_{k=1}^{j-1}c_{k}^{\left(j-1\right)}W_{k}\left(\boldsymbol{\Omega}\right),
\]
where the coefficients $c_{i}^{\left(j-1\right)}$ are uniquely chosen
to enforce the interpolation property
\[
W_{j}\left(\boldsymbol{\Omega}^{\left(i\right)}\right)=\sum_{k=1}^{j-1}c_{k}^{\left(j-1\right)}W_{k}\left(\boldsymbol{\Omega}^{\left(i\right)}\right),\,\,\,\,\,\,\,\,\,\,i=1,\ldots,j-1.
\]

To find the next interpolation point $\boldsymbol{\Omega}^{\left(j\right)}$,
we subtract from $W_{j}$ its interpolatory projection on to the previous
$j-1$ basis functions $W_{1},\ldots,W_{j-1}$:
\[
R_{j}\left(\boldsymbol{\Omega}\right)=W_{j}\left(\boldsymbol{\Omega}\right)-\left(P^{\left(j-1\right)}\mathbf{w}^{\left(j\right)}\right)\left(\boldsymbol{\Omega}\right).
\]
Then we define the next interpolation point via
\[
\boldsymbol{\Omega}_{j}=\text{argmax}_{\boldsymbol{\Omega}}\left|R_{j}\left(\boldsymbol{\Omega}\right)\right|.
\]
Note that, by construction,
\[
R_{j}\left(\boldsymbol{\Omega}^{\left(i\right)}\right)=W_{j}\left(\boldsymbol{\Omega}^{\left(i\right)}\right)-\sum_{i=1}^{j-1}c_{i}^{\left(j-1\right)}W_{i}\left(\boldsymbol{\Omega}^{\left(i\right)}\right)=0,\,\,\,\,\,\,\,\,\,\,i=1,\ldots,j-1.
\]
Therefore, $\boldsymbol{\Omega}^{\left(j\right)}\notin\left\{ \boldsymbol{\Omega}^{\left(1\right)},\ldots,\boldsymbol{\Omega}^{\left(j-1\right)}\right\} $.

\subsection{Collocation-based Dynamic Low Rank (DLR) Method}

The main algorithm starts with the initial representation
\[
\psi^{\left(0\right)}\left(\mathbf{x},\boldsymbol{\Omega}\right)=\psi\left(\mathbf{x},\boldsymbol{\Omega},t_{0}\right)=\boldsymbol{W}^{\text{T}}\left(\boldsymbol{\Omega},t_{0}\right)S\left(t_{0}\right)^{\text{T}}\mathbf{X}\left(\mathbf{x},t_{0}\right).
\]

We first select $r_{0}$ angles $\boldsymbol{\Omega}^{\left(j\right)}\left(t_{0}\right)$,
$j=1,\ldots,r{}_{0}$, using DEIM algorithm discussed in Section~\ref{sec:DEIM}
applied to the angular basis functions $\mathbf{W}\left(\boldsymbol{\Omega},t_{0}\right)$
at the beginning of the time step. We then evaluate the initial condition
for the next time step at these new angles,
\begin{align*}
\psi_{j}^{\left(0\right)}\left(\mathbf{x}\right) & =\psi\left(\mathbf{x},\boldsymbol{\Omega}^{\left(j\right)}\left(t_{0}\right),t_{0}\right)\\
 & =\sum_{k=1}^{r_{0}}L_{k}\left(\boldsymbol{\Omega}^{\left(j\right)}\left(t_{0}\right),t_{0}\right)X_{k}\left(\mathbf{x},t_{0}\right).
\end{align*}

Next, define 
\[
\psi_{j}\left(\mathbf{x},t\right)=\psi\left(\mathbf{x},\boldsymbol{\Omega}^{\left(j\right)}\left(t_{0}\right),t\right),\,\,\,\,\,\boldsymbol{\Psi}\left(\mathbf{x},t\right)=\left(\begin{array}{c}
\psi\left(\mathbf{x},\boldsymbol{\Omega}^{\left(1\right)}\left(t_{0}\right),t\right)\\
\vdots\\
\psi\left(\mathbf{x},\boldsymbol{\Omega}^{\left(r_{0}\right)}\left(t_{0}\right),t\right)
\end{array}\right).
\]
Then we solve the following equations (``column constraints'') for
$\boldsymbol{\Psi}$,
\begin{equation}
\boldsymbol{\Omega}^{\left(j\right)}\cdot\nabla_{\mathbf{x}}\psi_{j}+\left(\sigma+\frac{1}{c\Delta t}\right)\psi_{j}=\frac{1}{4\pi}\sigma_{s}\boldsymbol{\beta}^{\text{T}}\left(\hat{W}^{-1}\boldsymbol{\Psi}\right)+q+\frac{1}{c\Delta t}\psi_{j}^{\left(0\right)},\label{eq:steady state collocation}
\end{equation}
where we use the relationship derived in Section~\ref{subsec:Computing-the-scalar},
\[
\int_{\mathbb{S}^{2}}\psi\left(\mathbf{x},\boldsymbol{\Omega},t\right)d\boldsymbol{\Omega}\approx\boldsymbol{\beta}^{\text{T}}\left(\hat{W}^{-1}\boldsymbol{\Psi}\left(\mathbf{x},t\right)\right).
\]

In non-DLR production transport codes, the analog of equation (\ref{eq:steady state collocation}) is generally solved with ``source iteration'',
a type of fixed point iteration where the left hand side operator is inverted each iteration using fast sweep algorithms
while the entire right hand side is kept fixed for the inversion within an iteration but updated between them.  The sweep 
algorithms depend on angular information
being uncoupled and traveling along discrete angles on the left hand side, allowing the flow of dependencies to be traced in space
from the problem boundaries along each fixed direction, giving a fast triangular solve for each angle. This same structure holds in (\ref{eq:steady state collocation}),
allowing standard sweep and source iteration implementations to be carried over to the SN-like DLR method with minimal conversion.

Source iteration with sweeps along discrete angles has the physical interpretation of each iteration tracing light between scattering events,
with scattering occuring when the right hand side is updated between iterations. Therefore, in highly scattering media,
the number of iterations (scatter events) will be high before the final global distribution of light (for the time step) is fully converged
Diffusion synthetic acceleration, described in Section \ref{sec:DSA}, is used to accelerate convergence in the face of this.

After solving equation (\ref{eq:steady state collocation}), we orthonormalize
$\boldsymbol{\Psi}\left(\mathbf{x}\right)$ via stabilized Gram-Schmidt,
\[
\boldsymbol{\Psi}\left(\mathbf{x},t\right)=R\left(t\right)\mathbf{X}\left(\mathbf{x},t\right),
\]
where 
\[
\int_{\mathcal{X}}X_{i}\left(\mathbf{x}\right)X_{j}\left(\mathbf{x}\right)d\mathbf{x}=\delta_{i,j}.
\]

For ``row constraints'', we project equation (\ref{eq:linear transport with pseudoscattering})
using the new spatial basis functions,
\begin{equation}
\boldsymbol{\Omega}\cdot\left\langle \nabla_{\mathbf{x}}\psi,X_{i}\right\rangle _{\mathcal{X}}+\left\langle \left(\sigma+\frac{1}{c\Delta t}\right)\psi,X_{i}\right\rangle _{\mathcal{X}}=\frac{1}{4\pi}\left\langle \sigma_{s}\varphi,X_{i}\right\rangle _{\mathcal{X}}+\left\langle q+\frac{1}{c\Delta t}\psi^{\left(0\right)},X_{i}\right\rangle _{\mathcal{X}}.\label{eq:steady state-1}
\end{equation}
If we make the Galerkin approximation
\[
\psi\left(\mathbf{x},\boldsymbol{\Omega},t\right)=\sum_{j=1}^{r_{0}}L_{j}\left(\boldsymbol{\Omega},t\right)X_{j}\left(\mathbf{x},t\right),
\]
then we have the equations for $L_{j}\left(\boldsymbol{\Omega},t\right)$,
\begin{equation}
\left(A+M_{t}\right)\mathbf{L}=\frac{1}{4\pi}M_{s}\int_{\mathbb{S}^{2}}\mathbf{L}d\boldsymbol{\Omega}+\mathbf{Q},\label{eq:L step}
\end{equation}
where 
\[
\mathbf{L}\left(\boldsymbol{\Omega},t\right)=\left(\begin{array}{c}
L_{1}\left(\boldsymbol{\Omega},t\right)\\
\vdots\\
L_{r_{0}}\left(\boldsymbol{\Omega},t\right)
\end{array}\right),
\]
\[
A_{i,j}=\sum_{k}\Omega_{k}\left\langle \frac{\partial X_{j}}{\partial x_{k}},X_{i}\right\rangle _{\mathcal{X}},\,\,\,\,\,\,\left(M_{t}\right)_{i,j}=\left\langle \left(\sigma+\frac{1}{c\Delta t}\right)X_{j},X_{i}\right\rangle _{\mathcal{X}},
\]
and
\[
\left(M_{s}\right)_{i,j}=\left\langle \sigma_{s}X_{j},X_{i}\right\rangle _{\mathcal{X}},\,\,\,\,\,Q_{i}=\left\langle q+\frac{1}{c\Delta t}\psi^{\left(0\right)},X_{i}\right\rangle _{\mathcal{X}}.
\]
Note that we can compute $Q_{i}\left(\boldsymbol{\Omega},t\right)$
using the low-rank factorization of $\psi\left(\mathbf{x},\boldsymbol{\Omega},t_{0}\right)$:
\begin{align*}
Q_{i}\left(\boldsymbol{\Omega},t\right) & =\left\langle q\left(\cdot,t\right)+\frac{1}{c\Delta t}\psi\left(\cdot,\boldsymbol{\Omega},t_{0}\right),X_{i}\left(\cdot,t\right)\right\rangle _{\mathcal{X}}\\
 & =\left\langle q\left(\cdot,t\right),X_{i}\left(\cdot,t\right)\right\rangle _{\mathcal{X}}+\frac{1}{c\Delta t}\sum_{j=1}^{r_{0}}L_{i}\left(\boldsymbol{\Omega},t_{0}\right)\left\langle X_{i}\left(\cdot,t_{0}\right),X_{i}\left(\cdot,t\right)\right\rangle _{\mathcal{X}}.
\end{align*}

Next, orthonormalize the $L_{j}$ functions,
\[
L_{j}\left(\boldsymbol{\Omega},t\right)=\sum_{k=1}^{r_{0}}S_{jk}\left(t\right)W_{k}\left(\boldsymbol{\Omega},t\right)=\left(S\left(t\right)\boldsymbol{W}\left(t\right)\right)_{j}\left(\boldsymbol{\Omega},t\right),
\]
where 
\[
\int_{\mathbb{S}^{2}}W_{i}\left(\boldsymbol{\Omega}\right)W_{j}\left(\boldsymbol{\Omega}\right)d\boldsymbol{\Omega}=\delta_{i,j}.
\]
Then we finally obtain that
\begin{align*}
\psi\left(\mathbf{x},\boldsymbol{\Omega},t\right) & =\mathbf{L}^{\text{T}}\left(\boldsymbol{\Omega},t\right)\mathbf{X}\left(\mathbf{x},t\right)\\
 & =\boldsymbol{W}^{\text{T}}\left(\boldsymbol{\Omega},t\right)S\left(t\right)^{\text{T}}\mathbf{X}\left(\mathbf{x},t\right).
\end{align*}

We summarize the DLR solve in Algorithm~\ref{alg:SN-like-DLR-method}.
The input to this algorithm is initial condition,
\[
\psi^{\left(0\right)}\left(\mathbf{x},\boldsymbol{\Omega}\right)=\sum_{k=1}^{r_{0}}L_{k}^{\left(0\right)}\left(\boldsymbol{\Omega}\right)X_{k}^{\left(0\right)}\left(\mathbf{x}\right).
\]
In addition,
\[
L_{j}^{\left(0\right)}\left(\boldsymbol{\Omega}\right)=\sum_{k=1}^{r_{0}}S_{jk}^{\left(0\right)}W_{k}^{\left(0\right)}\left(\boldsymbol{\Omega}\right)=\left(S^{\left(0\right)}\boldsymbol{W}^{\left(0\right)}\right)_{j}\left(\boldsymbol{\Omega}\right),
\]
where 
\[
\int_{\mathbb{S}^{2}}W_{i}^{\left(0\right)}\left(\boldsymbol{\Omega}\right)W_{j}^{\left(0\right)}\left(\boldsymbol{\Omega}\right)d\boldsymbol{\Omega}=\delta_{i,j}.
\]

\begin{algorithm}
\caption{SN-like DLR method for the time step update\label{alg:SN-like-DLR-method}}

\begin{enumerate}
\item Update the effective scattering and transport source using equations
\ref{eq:effective scattering}-\ref{eq:effective source from Newton}
\item Select $r_{0}$ angles $\boldsymbol{\Omega}^{\left(j\right)}$, $j=1,\ldots,r{}_{0}$,
using DEIM algorithm in Section~\ref{sec:DEIM} applied to $\mathbf{W}^{\left(0\right)}\left(\boldsymbol{\Omega}\right)$.
\item Evaluate the initial condition at the new angles, for $j=1,\ldots,r{}_{0}$,
\[
\psi_{j}^{\left(0\right)}\left(\mathbf{x}\right)=\psi^{\left(0\right)}\left(\mathbf{x},\boldsymbol{\Omega}^{\left(j\right)}\right)=\sum_{k=1}^{r_{0}}L_{k}^{\left(0\right)}\left(\boldsymbol{\Omega}^{\left(j\right)}\right)X_{k}^{\left(0\right)}\left(\mathbf{x}\right).
\]
\item Using fixed-point iteration, solve the SN-like system of $r_{0}$ equations
\[
\boldsymbol{\Omega}^{\left(j\right)}\cdot\nabla_{\mathbf{x}}\psi_{j}+\left(\sigma+\frac{1}{c\Delta t}\right)\psi_{j}=\frac{1}{4\pi}\sigma_{s}\boldsymbol{\beta}^{\text{T}}\left(\hat{W}^{-1}\boldsymbol{\Psi}\right)+q+\frac{1}{c\Delta t}\psi_{j}^{\left(0\right)},
\]
where $\hat{W}_{ij}=W_{j}\left(\boldsymbol{\Omega}^{\left(i\right)}\right)$,
$\beta_{j}=\int_{\mathbb{S}^{2}}W_{j}\left(\boldsymbol{\Omega}\right)d\boldsymbol{\Omega}$,
\[
\psi_{j}\left(\mathbf{x}\right)=\psi\left(\mathbf{x},\boldsymbol{\Omega}^{\left(j\right)}\right),\,\,\,\,\,\boldsymbol{\Psi}\left(\mathbf{x}\right)=\left(\begin{array}{c}
\psi_{1}\left(\mathbf{x}\right)\\
\vdots\\
\psi_{r_{0}}\left(\mathbf{x}\right)
\end{array}\right).
\]
Sweep algorithms are used to invert the left hand side operator each iteration.
Applying DSA from Section \ref{sec:DSA} keeps the iteration count low.

\item Orthonormalize $\boldsymbol{\Psi}\left(\mathbf{x}\right)$ via stabilized
Gram-Schmidt,
\[
\boldsymbol{\Psi}\left(\mathbf{x}\right)=R\mathbf{X}\left(\mathbf{x}\right),\,\,\,\,\,\,\,\,\,\,\,\,\int_{\mathcal{X}}X_{i}\left(\mathbf{x}\right)X_{j}\left(\mathbf{x}\right)d\mathbf{x}=\delta_{i,j}.
\]
\item Project the initial condition $\psi_{j}^{\left(0\right)}$ on to the
$r$ spatial basis functions $X_{j}$,
\[
\tilde{L}_{j}^{\left(0\right)}=\left\langle \psi^{\left(0\right)},X_{j}\right\rangle _{\mathcal{X}}=\sum_{j=1}^{r_{0}}L_{i}^{\left(0\right)}\left\langle X_{i}^{\left(0\right)},X_{j}\right\rangle _{\mathcal{X}}.
\]
\item Solve 
\[
\left(A+M_{t}\right)\mathbf{L}=\frac{1}{4\pi}M_{s}\int_{\mathbb{S}^{2}}\mathbf{L}d\boldsymbol{\Omega}+\mathbf{Q},
\]
where
\[
A_{i,j}=\sum_{k}\Omega_{k}\left\langle \frac{\partial X_{j}}{\partial x_{k}},X_{i}\right\rangle _{\mathcal{X}},\,\,\,\,\,\,\left(M_{t}\right)_{i,j}=\left\langle \left(\sigma+\frac{1}{c\Delta t}\right)X_{j},X_{i}\right\rangle _{\mathcal{X}},
\]
\[
\left(M_{s}\right)_{i,j}=\left\langle \sigma_{s}X_{j},X_{i}\right\rangle _{\mathcal{X}},\,\,\,\,\,Q_{i}=\left\langle q,X_{i}\right\rangle _{\mathcal{X}}+\frac{1}{c\Delta t}\tilde{L}_{i}^{\left(0\right)}.
\]
\item Orthonormalize the $L_{j}$ functions, $j=1,\ldots,r$,
\[
L_{j}\left(\boldsymbol{\Omega}\right)=\sum_{k=1}^{r_{0}}S_{jk}W_{k}\left(\boldsymbol{\Omega}\right)=\left(S\boldsymbol{W}\right)_{j}\left(\boldsymbol{\Omega}\right),\,\,\,\,\,\,\,\,\,\,\,\,,\int_{\mathbb{S}^{2}}W_{i}W_{j}d\boldsymbol{\Omega}=\delta_{i,j}.
\]
\item Update the material temperature using equation (\ref{eq:temperature update})
and the new mean intensity $\varphi$
\end{enumerate}
\end{algorithm}

\subsection{Spatial discretization for collocation solve} \label{subsec:spatial discretization}

In preparation for defining the DSA algorithm in the next section to accelerate the solve for (\ref{eq:steady state collocation}),
we here detail the spacial discretization of the equations, as it is imperative that DSA be derived consistently with it (see \citep{ADAMS20023}).

Define
\begin{equation}
\mathbf{v}^{\text{T}}G_{j}\mathbf{u}=-\sum_{e}\int_{\kappa_{e}}\frac{\partial v}{\partial x_{j}}ud\mathbf{x},\,\,\,\,j=1,2,3.\label{eq:G_j}
\end{equation}
Also, define the angle-dependent face matrix
\[
\mathbf{v}^{\text{T}}F_{d}\mathbf{u}=-\sum_{f}\int_{\Gamma_{f}}\boldsymbol{\Omega}^{\left(d\right)}\cdot\mathbf{n}\left[\left[v\right]\right]\left\{ \left\{ u\right\} \right\} dS+\frac{1}{2}\sum_{f}\int_{\Gamma_{f}}\left|\boldsymbol{\Omega}^{\left(d\right)}\cdot\mathbf{n}\right|\left[\left[u\right]\right]\left[\left[v\right]\right]dS.
\]
Finally, define the following mass matrices 
\[
\mathbf{v}^{\text{T}}M_{t}\mathbf{u}=-\sum_{e}\int_{\kappa_{e}}\sigma_{t}uvd\mathbf{x},
\]
\[
\mathbf{v}^{\text{T}}M_{a}\mathbf{u}=-\sum_{e}\int_{\kappa_{e}}\sigma_{a}uvd\mathbf{x},
\]
\[
\mathbf{v}^{\text{T}}M\mathbf{u}=-\sum_{e}\int_{\kappa_{e}}uvd\mathbf{x}.
\]
Then equation (\ref{eq:steady state collocation}) is discretized
by
\[
\left(\boldsymbol{\Omega}^{\left(d\right)}\cdot\mathbf{G}+F_{d}\right)\boldsymbol{\psi}_{d}+\left(M_{a}+\frac{1}{c\Delta t}M\right)\boldsymbol{\psi}_{d}=\frac{1}{4\pi}M_{s}\boldsymbol{\beta}^{\text{T}}\left(\hat{W}^{-1}\boldsymbol{\psi}\right)+\mathbf{q}_{d}+\frac{1}{c\Delta t}M\boldsymbol{\psi}_{0,d}.
\]
Here 
\[
\boldsymbol{\psi}=\left(\begin{array}{c}
\boldsymbol{\psi}_{1}\\
\vdots\\
\boldsymbol{\psi}_{r_{0}}
\end{array}\right).
\]

\subsection{Diffusion Synthetic Acceleration (DSA) for SN-like DLR method} \label{sec:DSA}

It is well-understood that, for ``optically thick'' problems that
are characterized by small photon mean free path relative to the mesh
resolution, a basic fixed-point iteration to solve spatially discretized
version of equation (\ref{eq:linear transport with pseudoscattering})
converges arbitrarily slowly (in fact, like the inverse square of
the photon mean free path). Diffusion Synthetic Acceleration (DSA)
(cf. \citep{Kopp1963}, \citep{LEBEDEV-1969}, and \citep{Larsen1984})
is the traditional means to restore fast convergence; see \citep{ADAMS20023}
for an excellent discussion of iterative methods for solving linear
transport equations such as (\ref{eq:linear transport with pseudoscattering}),
and how DSA gives convergence independent of the photon mean free
path.

The need for DSA acceleration carries over to the equation (\ref{eq:steady state collocation}).
It is crucial for iterative efficiency (see \citep{ADAMS20023}) that
a so-called consistent DSA method is derived. We briefly outline how
DSA can be used, with little modification, to efficiently accelerate
the fixed-point solution to equation (\ref{eq:steady state collocation});
a detailed derivation is included in Section~\ref{subsec:Derivation-of-DSA}
of the appendix.

Introduce the standard diffusion limit scaling in equation (\ref{eq:steady state collocation})
using the non-dimensional parameter $0<\varepsilon\leq1$,

\begin{equation}
\left(\boldsymbol{\Omega}^{\left(d\right)}\cdot\mathbf{G}+F_{d}\right)\boldsymbol{\psi}_{d}+\varepsilon^{-1}M_{t}\boldsymbol{\psi}_{d}=\frac{1}{4\pi}\left(\varepsilon^{-1}M_{t}-\varepsilon M_{a}\right)\boldsymbol{\phi}+\varepsilon MQ_{d},\label{eq:steady state collocation, diffusion scaling}
\end{equation}
where
\[
\boldsymbol{\phi}=\boldsymbol{\beta}^{\text{T}}\left(\hat{W}^{-1}\boldsymbol{\psi}\right).
\]
Here $\varepsilon$ characterizes the photon mean free path relative
to the mesh spacing, and $0<\varepsilon\ll1$ results in arbitrarily
large iteration counts for fixed point iteration of equation (\ref{eq:steady state collocation, diffusion scaling}).

In Section~\ref{subsec:Derivation-of-DSA}, we derive the following
discrete diffusion approximation,

\[
\frac{1}{3}\sum_{j=1}^{3}\left(G_{j}+\overline{F}_{j}\right)^{\text{T}}M_{t}^{-1}\left(G_{j}+\overline{F}_{j}\right)\boldsymbol{\phi}+F_{0}\boldsymbol{\phi}+M_{a}\boldsymbol{\phi}\approx MQ^{+}.
\]
Here $G_{j}$ is defined in equation (\ref{eq:G_j}) and $\overline{F}_{j}$
is defined by
\[
\mathbf{v}^{\text{T}}\overline{F}_{j}\mathbf{u}=\sum_{\Gamma\in\mathcal{F}}\int_{\Gamma}\left\llbracket u\right\rrbracket \left\{ v\right\} dS,
\]
and 
\begin{equation}
\mathbf{v}^{\text{T}}F_{0}\mathbf{u}=\frac{1}{4\pi}\sum_{\Gamma\in\mathcal{F}}\left(\frac{1}{2}\int_{\Gamma}\boldsymbol{\beta}^{\text{T}}\hat{W}^{-1}\left[\left|\boldsymbol{\Omega}^{\left(d\right)}\cdot\mathbf{n}\right|\right]_{d=1}^{r}\right)\left\llbracket u\right\rrbracket \left\llbracket v\right\rrbracket dS.\label{eq:penalty term}
\end{equation}
We note that this form is identical to a traditional consistent DSA
discretization for DG discretizations of equation (\ref{eq:linear transport with pseudoscattering}).
The only difference is in the penalty term (\ref{eq:penalty term})
(cf. \citep{Haut2020}); in practice, we have observed that using
the penalty term from \citep{Haut2020} is sufficient for good acceleration,
and we use this version for simplicity in the numerical experiments
in Section~\ref{sec:Numerical-experiments}.

\section{Numerical experiments\label{sec:Numerical-experiments}}

Here we compare the PN-like and SN-like DLR methods against a classic SN TRT solver on two challenging test problems.
The test problem in Section \ref{subsec:lattice} contains extremely optically thick and optically thin spatial regions, and its efficient solution necessitates implicit time-stepping and---in the SN-like DLR method---carefully designed acceleration methods such as the DSA scheme discussed in Section~\ref{sec:DSA}. The test problem in Section \ref{subsec:Hohlraum problem} is even more challenging for standard SN methods, due to the point-like radiation sources from the hot spots on the hohlraum wall in to the very thin hohlraum cavity, and the even more severe ``ray effects" resulting from the under-resolution in angle. 

The spatial discretization we used for the SN-like DLR method is detailed in Section~\ref{subsec:spatial discretization}. For the spatial discretization of the PN-like DLR method, we use a second-order continuous Galerkin finite element discretization.  All three methods (the SN, SN-like DLR, and PN-like DLR) use the same spatial meshes, which are displayed in Figures~\ref{fig:lattice_materials}~and~\ref{fig:hohlraum_problem}.

Comparing the relative efficiency of the PN-like and SN-like methods is beyond the scope of this paper. However, we do make the following observations. First, we observe that the SN-like DLR method has the same computational performance as the SN method when both use the same number of angles; this results from both using the same highly efficient sweep algorithms and both requiring nearly identical numbers of iterations to converge the implicit solution for each time step. Importantly, we also observe that the SN-like DLR method displays noticeably less angular artifacts from under-resolution than the SN method when both use a comparable number of angles. In contrast, the PN-like DLR method exhibits no ray effects; although it requires two $r_0 \times r_0$ block-diffusion system solves per time step, this cost is likely offset in problems where ray effects are particularly challenging, such as the test problem in Section~\ref{subsec:Hohlraum problem}. We leave the careful comparison of the relative efficiency of these schemes to future work.

\subsection{Lattice problem} \label{subsec:lattice}

We consider a grey version of the problem in \citep{brunner2023family}.
This problem has highly heterogeneous materials---some very optically
thick and some very optically thin---and serves as an excellent stress
test for numerical TRT methods.

As shown in Figure~\ref{fig:lattice_materials}, the lattice
problem consists of a checkerboard pattern of materials with very
optically thick iron blocks (the light blue-green region corresponds
to hot iron and the yellow regions to cold iron), moderately thick
blocks of diamond (red regions), and optically thin blocks of foam
(blue regions). The mesh elements are chosen to accumulate at material
interfaces (see Figure~\ref{fig:lattice_materials}).

For simplicity, the opacity of each material is chosen to be fixed
instead of temperature dependent; however, realistic opacity values
are chosen from the opacity plots in figure $3$ in \citep{brunner2023family}.
Specifically, the blocks of iron (light blue-green and yellow blocks)
have material opacity $\sigma_{a}=10^{4}\text{\text{cm}}^{-1}$, specific
heat $c_{v}=0.05427\,\text{GJ}/(\text{g}\cdot\text{keV})$, and density
$\rho=8.0\text{g}/\text{\text{cm}}^{3}$. The blocks of diamond (red
blocks) have $\sigma_{a}=10^{2}\text{\text{cm}}^{-1}$, $c_{v}=0.05427\,\text{GJ}/(\text{g}\cdot\text{keV})$,
and $\rho=2.0\text{g}/\text{\text{cm}}^{3}$. The blocks of foam (blue
blocks) have $\sigma_{a}=4\text{\text{cm}}^{-1}$, $c_{v}=0.02412\,\text{GJ}/(\text{g}\cdot\text{keV})$,
and $\rho=0.2\text{g}/\text{\text{cm}}^{3}$.

Figure~\ref{fig:lattice_dlr_comparison} shows the PN-like
DLR scheme from Section~\ref{sec:PN like DLR scheme} with rank $8$,
the SN-like DLR scheme from Section~\ref{sec:SN like DLR scheme}
with rank $16$, and the S6 scheme (with $36$ angles) and S10 scheme (with $100$ angles) applied to
the lattice problem. Plots Figure~\ref{fig:lattice_dlr_comparison}(a)-(d)
show the material temperature at times $t=3$, $t=14$, and $t=25$.
As can be seen from the plots, the S6 method (Figure~\ref{fig:lattice_dlr_comparison}(c),(g),(k))
has much more visible ``ray effects'' from under-resolution in angle
than either of the DLR methods. Note that the SN-like DLR method with
$16$ angles has about the computational cost and memory footprint
of the S6 method (by symmetry in xy geometry, we need only solving for half the angles for S6 and S10).

\begin{figure}[htbp]
\centering
\includegraphics[width=0.35\textwidth]{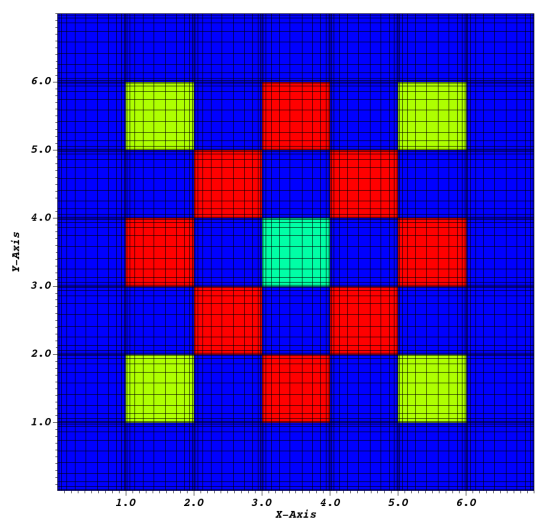}
\caption{
Mesh and materials for the lattice problem defined in \citep{brunner2023family}.
The blue-green block in the center of the spatial domain corresponds to hot iron;
the red blocks correspond to diamond; the yellow blocks correspond to cold iron;
and the blue blocks correspond to foam.
}
\label{fig:lattice_materials}
\end{figure}

\begin{figure}[htbp]
\centering
\setlength{\tabcolsep}{2pt}

\begin{minipage}[t]{0.23\textwidth}\centering
\includegraphics[width=\linewidth]{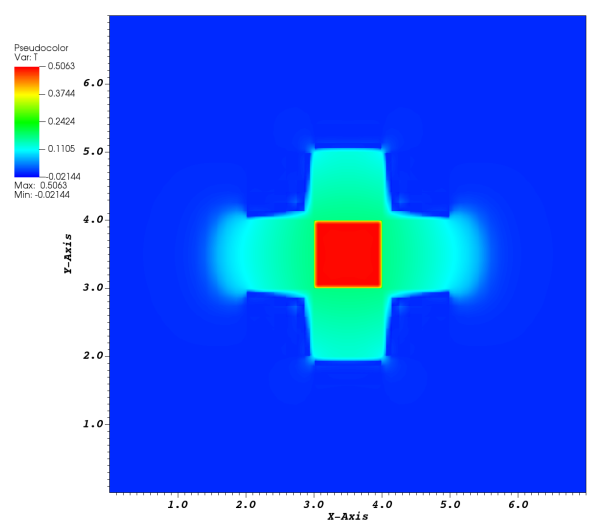}\par
\vspace{1pt}
{\footnotesize (a) PN-like DLR, rank $8$}
\end{minipage}\hfill
\begin{minipage}[t]{0.23\textwidth}\centering
\includegraphics[width=\linewidth]{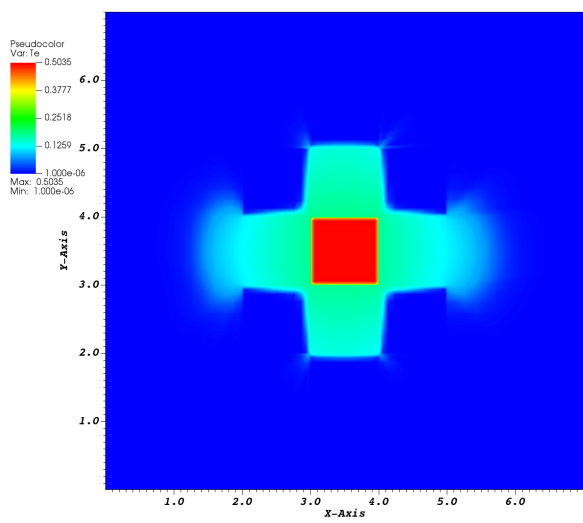}\par
\vspace{1pt}
{\footnotesize (b) SN-like DLR, rank $16$}
\end{minipage}\hfill
\begin{minipage}[t]{0.23\textwidth}\centering
\includegraphics[width=\linewidth]{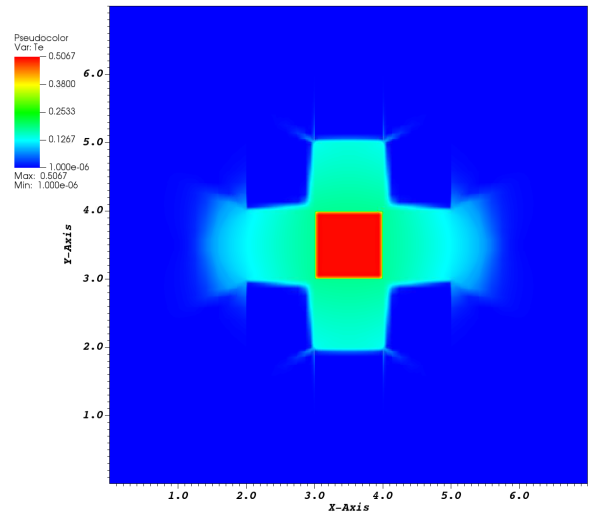}\par
\vspace{1pt}
{\footnotesize (c) S6 ($18$ angles)}
\end{minipage}\hfill
\begin{minipage}[t]{0.23\textwidth}\centering
\includegraphics[width=\linewidth]{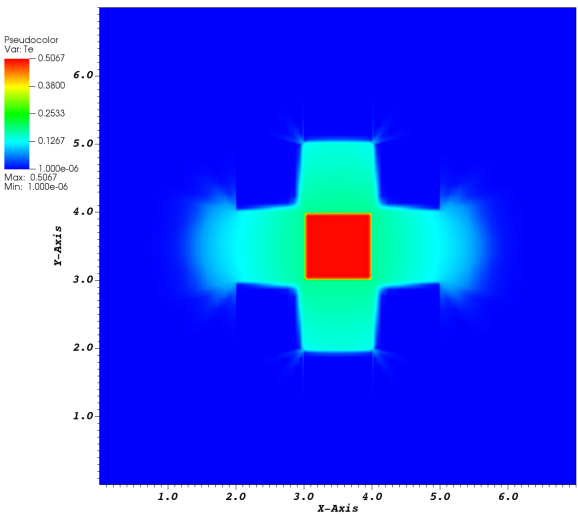}\par
\vspace{1pt}
{\footnotesize (d) S10 ($50$ angles)}
\end{minipage}

\vspace{3mm}

\begin{minipage}[t]{0.23\textwidth}\centering
\includegraphics[width=\linewidth]{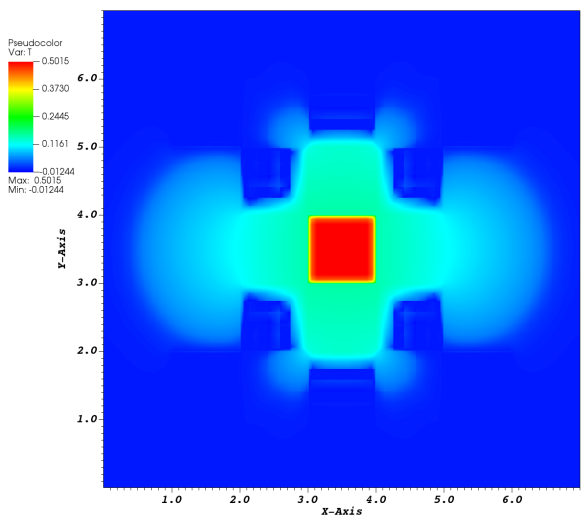}\par
\vspace{1pt}
{\footnotesize (e) PN-like DLR, rank $8$}
\end{minipage}\hfill
\begin{minipage}[t]{0.23\textwidth}\centering
\includegraphics[width=\linewidth]{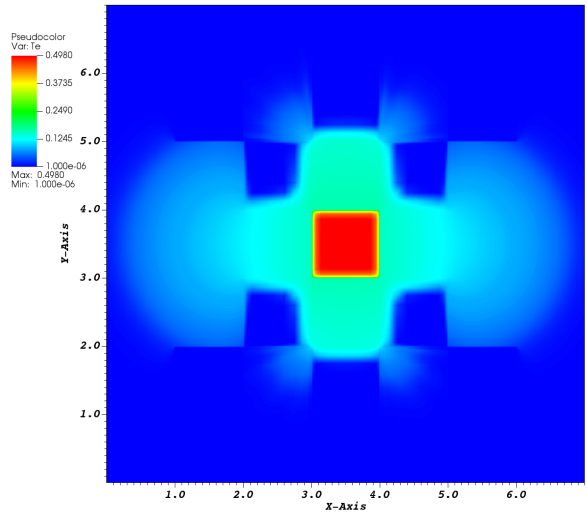}\par
\vspace{1pt}
{\footnotesize (f) SN-like DLR, rank $16$}
\end{minipage}\hfill
\begin{minipage}[t]{0.23\textwidth}\centering
\includegraphics[width=\linewidth]{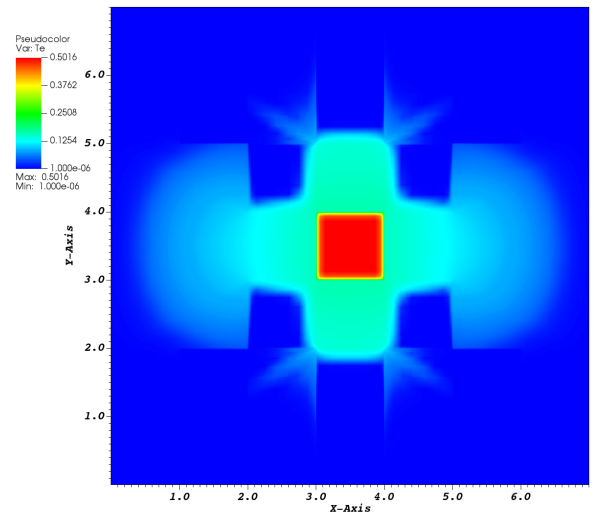}\par
\vspace{1pt}
{\footnotesize (g) S6 ($18$ angles)}
\end{minipage}\hfill
\begin{minipage}[t]{0.23\textwidth}\centering
\includegraphics[width=\linewidth]{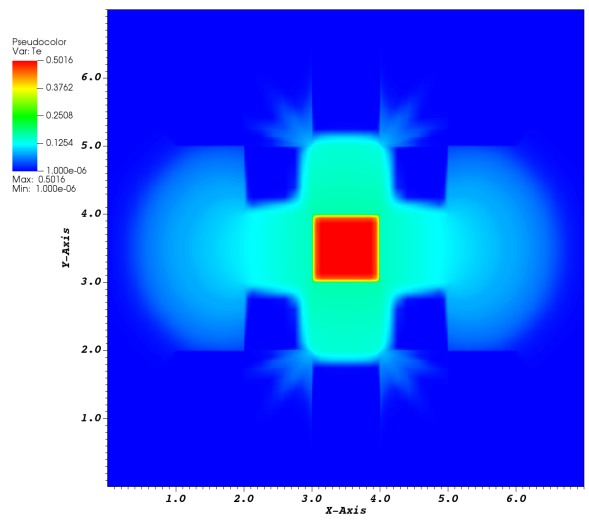}\par
\vspace{1pt}
{\footnotesize (h) S10 ($50$ angles)}
\end{minipage}

\vspace{3mm}

\begin{minipage}[t]{0.23\textwidth}\centering
\includegraphics[width=\linewidth]{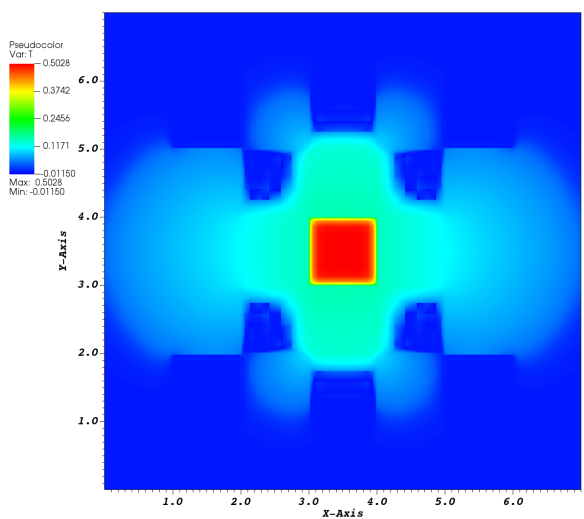}\par
\vspace{1pt}
{\footnotesize (i) PN-like DLR, rank $8$}
\end{minipage}\hfill
\begin{minipage}[t]{0.23\textwidth}\centering
\includegraphics[width=\linewidth]{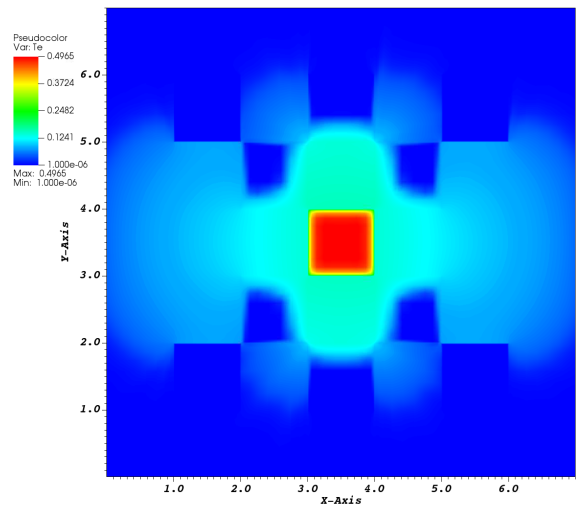}\par
\vspace{1pt}
{\footnotesize (j) SN-like DLR, rank $16$}
\end{minipage}\hfill
\begin{minipage}[t]{0.23\textwidth}\centering
\includegraphics[width=\linewidth]{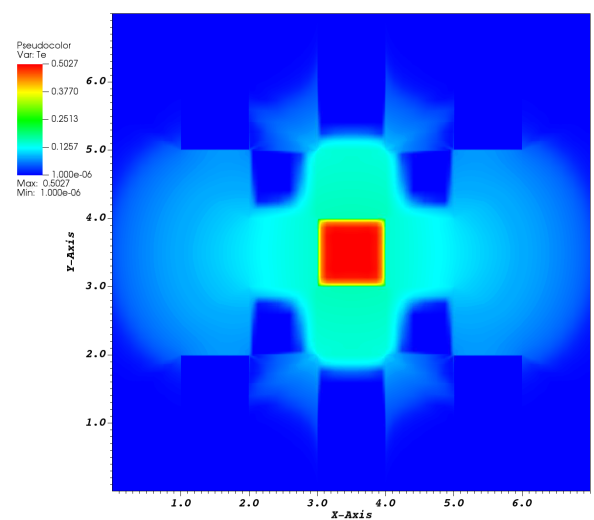}\par
\vspace{1pt}
{\footnotesize (k) S6 ($18$ angles)}
\end{minipage}\hfill
\begin{minipage}[t]{0.23\textwidth}\centering
\includegraphics[width=\linewidth]{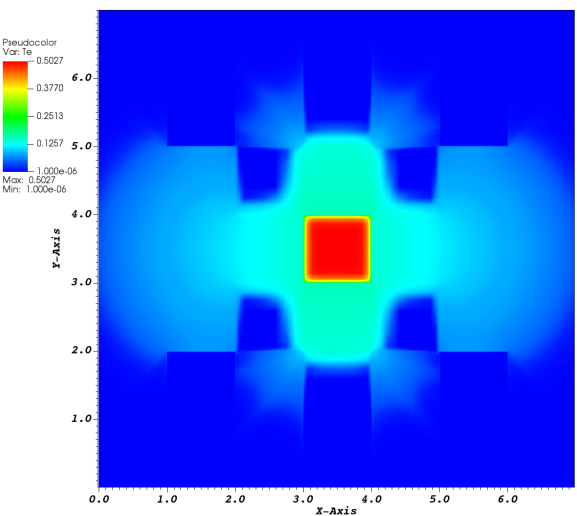}\par
\vspace{1pt}
{\footnotesize (l) S10 ($50$ angles)}
\end{minipage}

\caption{%
PN-like DLR method with rank $8$, SN-like DLR method with rank $16$,
S6 ($18$ angles), and S10 ($50$ angles) for the lattice problem.
Material temperature at times $t=3$, $t=14$, and $t=25$.
}
\label{fig:lattice_dlr_comparison}
\end{figure}

\subsection{Hohlraum problem\label{subsec:Hohlraum problem}}

We consider a hohlraum problem, as shown in Figure~\ref{fig:hohlraum_problem}.
Figure~\ref{fig:hohlraum_problem}(a) shows the material regions
for the hohlraum problem. The center light blue disc in Figure~\ref{fig:hohlraum_problem}(a)
corresponds to the capsule, and has an absorption opacity of $\sigma_{a}=10^{2}\,\text{cm}^{-1}$,
a specific heat of $c_{v}=0.05427\text{\,GJ}/(\text{g}\cdot\text{keV})$,
and a density of $\rho=2.0\,\text{g}/\text{\text{cm}}^{3}$. The dark
blue region surrounding the capsule corresponds to a foam-filled hohlraum
with an absorption opacity of $\sigma_{a}=4\,\text{cm}^{-1}$, a specific
heat of $c_{v}=0.02412\,\text{GJ}/(\text{g}\cdot\text{keV})$, and
a density of $\rho=0.2\,\text{g}/\text{\text{cm}}^{3}$. The red region
in Figure~\ref{fig:hohlraum_problem}(a) corresponds to the hohlraum
walls and has an absorption opacity of $\sigma_{a}=10^{4}\,\text{cm}^{-1}$,
a specific heat of $c_{v}=0.05427\,\text{GJ}/(\text{g}\cdot\text{keV})$,
and a density of $\rho=8.0\text{g}/\text{\text{cm}}^{3}$. As shown
in Figure~\ref{fig:hohlraum_problem}(b), the initial material temperature
is $T_{\text{init}}=0.0001\,\text{keV}$ in the blue regions and $T_{\text{init}}=0.3\,\text{keV}$
in the eight ``hot spots'' in the hohlraum walls (shown in the red
regions). The hot spots represent where the cones of laser beams entering
the hohlraum through the wall openings heat the outer surface of the
hohlraum walls.

Figure~\ref{fig:hohlraum_dlr} shows the PN-like DLR scheme from Section~\ref{sec:PN like DLR scheme}
with rank $8$, the SN-like DLR scheme from Section~\ref{sec:SN like DLR scheme}
with rank $16$, and the S6 scheme (with $36$ angles) applied to
the hohlraum problem. Because the point-source-like hot spots on the
hohlraum walls radiate into the very optically thin cavity, the ray
effects in both the SN-like DLR method and the standard SN method
are highly visible; however, the ray effects in the SN-like DLR method
appear somewhat mitigated in comparison to the S6 method. In contrast,
the PN-like DLR method shows a much smoother behavior in angle.

\begin{figure}[htbp]
\centering

\begin{minipage}[t]{0.30\textwidth}\centering
\includegraphics[width=\linewidth]{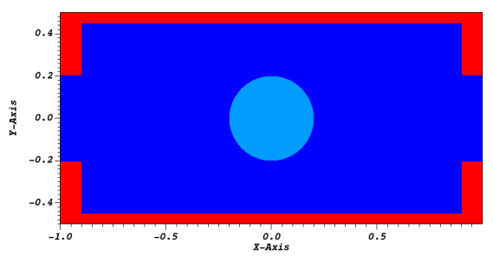}\par
\vspace{1pt}
{\footnotesize (a) Hohlraum materials}
\end{minipage}\hfill
\begin{minipage}[t]{0.30\textwidth}\centering
\includegraphics[width=\linewidth]{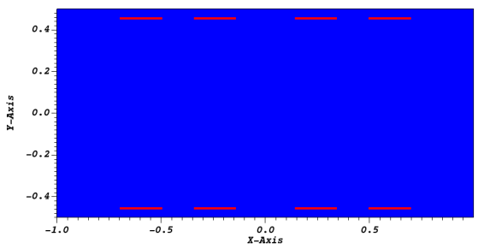}\par
\vspace{1pt}
{\footnotesize (b) Initial material temperature}
\end{minipage}\hfill
\begin{minipage}[t]{0.30\textwidth}\centering
\includegraphics[width=\linewidth]{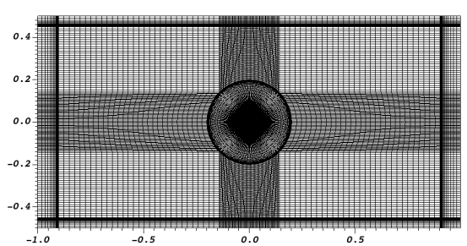}\par
\vspace{1pt}
{\footnotesize (c) Spatial mesh}
\end{minipage}

\caption{
Material regions, mesh, and initial conditions for the hohlraum problem.
}
\label{fig:hohlraum_problem}
\end{figure}

\begin{figure}[htbp]
\centering

\begin{minipage}[t]{0.30\textwidth}\centering
\includegraphics[width=\linewidth]{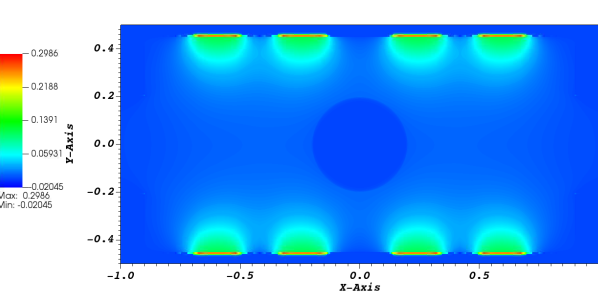}\par
\vspace{1pt}
{\footnotesize (a) PN-like DLR, rank $8$}
\end{minipage}\hfill
\begin{minipage}[t]{0.30\textwidth}\centering
\includegraphics[width=\linewidth]{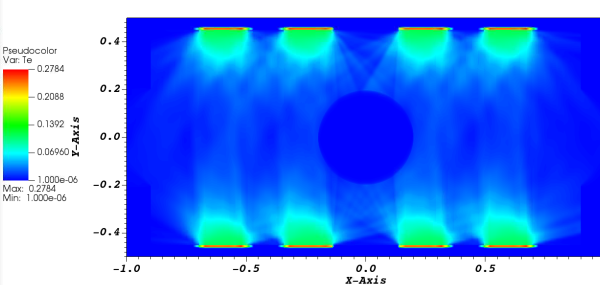}\par
\vspace{1pt}
{\footnotesize (b) SN-like DLR, rank $16$}
\end{minipage}\hfill
\begin{minipage}[t]{0.30\textwidth}\centering
\includegraphics[width=\linewidth]{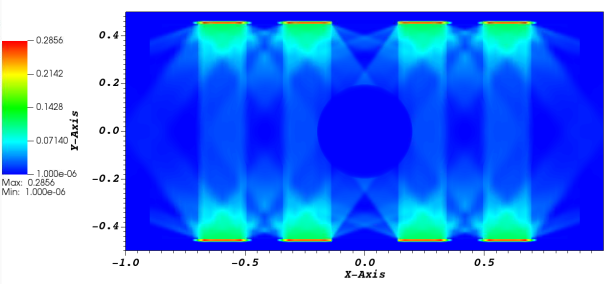}\par
\vspace{1pt}
{\footnotesize (c) S6 ($36$ angles)}
\end{minipage}

\vspace{3mm}

\begin{minipage}[t]{0.30\textwidth}\centering
\includegraphics[width=\linewidth]{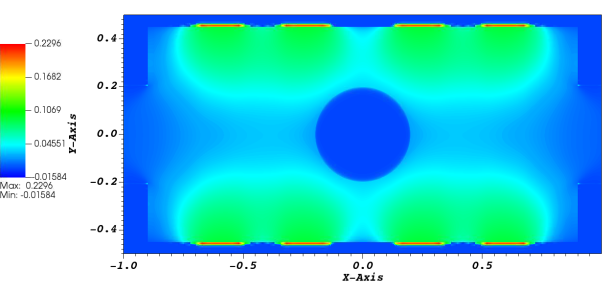}\par
\vspace{1pt}
{\footnotesize (d) PN-like DLR, rank $8$}
\end{minipage}\hfill
\begin{minipage}[t]{0.30\textwidth}\centering
\includegraphics[width=\linewidth]{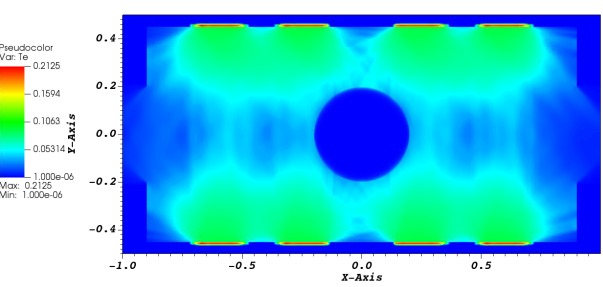}\par
\vspace{1pt}
{\footnotesize (e) SN-like DLR, rank $16$}
\end{minipage}\hfill
\begin{minipage}[t]{0.30\textwidth}\centering
\includegraphics[width=\linewidth]{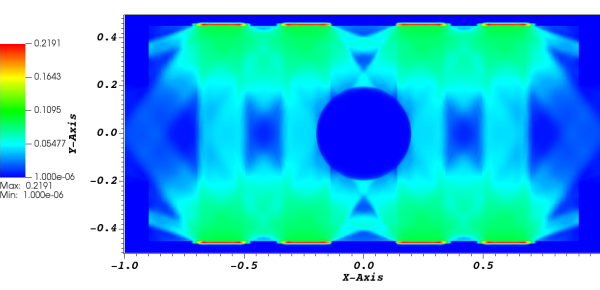}\par
\vspace{1pt}
{\footnotesize (f) S6 ($36$ angles)}
\end{minipage}

\vspace{3mm}

\begin{minipage}[t]{0.30\textwidth}\centering
\includegraphics[width=\linewidth]{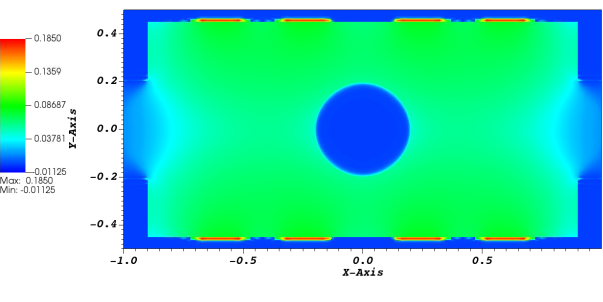}\par
\vspace{1pt}
{\footnotesize (g) PN-like DLR, rank $8$}
\end{minipage}\hfill
\begin{minipage}[t]{0.30\textwidth}\centering
\includegraphics[width=\linewidth]{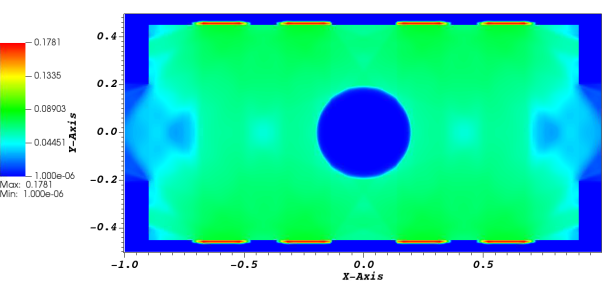}\par
\vspace{1pt}
{\footnotesize (h) SN-like DLR, rank $16$}
\end{minipage}\hfill
\begin{minipage}[t]{0.30\textwidth}\centering
\includegraphics[width=\linewidth]{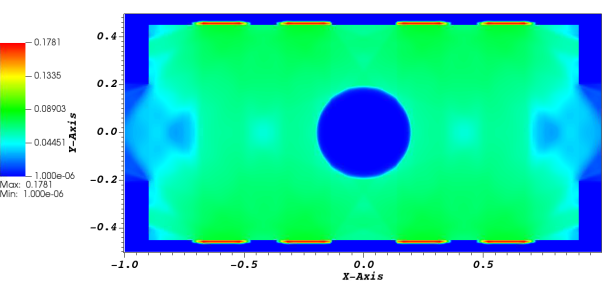}\par
\vspace{1pt}
{\footnotesize (i) S6 ($36$ angles)}
\end{minipage}

\caption{
PN-like DLR method with rank $8$, SN-like DLR method with rank $16$,
and S6 ($36$ angles) for the hohlraum problem.
Plots show the material temperature at times $t=0.3$, $t=1.2$, and $t=3.0$.
}
\label{fig:hohlraum_dlr}
\end{figure}

\section{Conclusions}
We develop two new types of DLR schemes for TRT.  One scheme takes a traditional DLR approach applied to the TRT equations,
resulting in a PN-like scheme.  The novelty of this method is that it uses an even-parity formulation to allow the time stepping
of the method to be done efficiently in the lack of traditional fast transport sweeps.  The second method uses a collocation-DLR approach
to produce an SN-like DLR scheme.  This method has the advantage that it can be adapted to current SN production methods that use
state-of-the-art sweep solvers with minimal conversion.  In combination with DLR compression, this approach is particularly promising
for adoption in application codes.

We demonstrate the two new methods in comparison with a traditional non-DLR TRT SN solver on two challenging test problems.  The results show
that the methods give solutions of comparable or superior quality to the non-DLR method even at modest rank.

\section{Appendix}

\subsection{Derivation of time-stepping scheme for PN-like DLR method\label{subsec:Appendix A}}

Assume that 
\[
\frac{\partial\psi}{\partial t}=\mathcal{L}\psi.
\]
Also assume that
\begin{align*}
\psi\left(\mathbf{x},\boldsymbol{\Omega},t\right) & =\sum_{i,j=1}^{r}X_{i}\left(\mathbf{x},t\right)S_{ij}\left(t\right)W_{j}\left(\boldsymbol{\Omega},t\right)\\
 & =\sum_{j=1}^{r}K_{j}\left(\mathbf{x},t\right)W_{j}\left(\boldsymbol{\Omega},t\right)\\
 & =\sum_{i=1}^{r}X_{i}\left(\mathbf{x},t\right)L_{i}\left(\boldsymbol{\Omega},t\right),
\end{align*}
where
\[
K_{j}=\sum_{i=1}^{r}X_{i}S_{ij},\,\,\,\,\,L_{i}=\sum_{j=1}^{r}S_{ij}W_{j}.
\]

Now
\begin{align*}
\frac{\partial\psi}{\partial t} & =\sum_{j=1}^{r}\dot{K}_{j}W_{j}+\sum_{j=1}^{r}K_{j}\dot{W}_{j}\\
 & =\sum_{j=1}^{r}\dot{X}_{j}L_{j}+\sum_{j=1}^{r}X_{j}\dot{L}_{j}.
\end{align*}
Therefore,
\[
\dot{L}_{j}=\left\langle \frac{\partial\psi}{\partial t},X_{j}\right\rangle _{\mathcal{X}},\,\,\,\,\,\dot{K}_{j}=\left\langle \frac{\partial\psi}{\partial t},W_{j}\right\rangle _{\mathbb{S}^{2}}.
\]
So for consistency
\[
\dot{L}_{j}=\left\langle \frac{\partial\psi}{\partial t},X_{j}\right\rangle _{\mathcal{X}}=\left\langle \mathcal{L}\psi,X_{j}\right\rangle _{\mathcal{X}},
\]
and
\[
\dot{K}_{j}=\left\langle \frac{\partial\psi}{\partial t},W_{j}\right\rangle _{\mathcal{X}}=\left\langle \mathcal{L}\psi,W_{j}\right\rangle _{\mathcal{X}}.
\]
That is, we compute $\dot{K}_{j}$, $\dot{L}_{j}$ from Galerkin conditions
\[
\left\langle \frac{\partial\psi}{\partial t}-\mathcal{L}\psi,X_{j}\right\rangle _{\mathcal{X}}=0,\,\,\,\,\,\left\langle \frac{\partial\psi}{\partial t}-\mathcal{L}\psi,W_{j}\right\rangle _{\mathbb{S}^{2}}=0.
\]
Now,
\[
\left\langle \mathcal{L}\psi,X_{j}\right\rangle _{\mathcal{X}}=\sum_{i=1}^{r}\left\langle \mathcal{L}\left(X_{i}L_{i}\right),X_{j}\right\rangle _{\mathcal{X}},
\]
\[
\left\langle \mathcal{L}\psi,W_{j}\right\rangle _{\mathcal{X}}=\sum_{i=1}^{r}\left\langle \mathcal{L}\left(W_{i}K_{i}\right),W_{j}\right\rangle _{\mathbb{S}^{2}},
\]
and so
\[
\dot{L}_{j}=\sum_{i=1}^{r}\left\langle \mathcal{L}\left(X_{i}L_{i}\right),X_{j}\right\rangle _{\mathcal{X}},
\]
\[
\dot{K}_{j}=\sum_{i=1}^{r}\left\langle \mathcal{L}\left(W_{i}K_{i}\right),W_{j}\right\rangle _{\mathbb{S}^{2}}.
\]

Finally consider
\[
\psi\left(\mathbf{x},\boldsymbol{\Omega},t\right)=\sum_{i,j=1}^{r}X_{i}\left(\mathbf{x},t\right)S_{ij}\left(t\right)W_{j}\left(\boldsymbol{\Omega},t\right).
\]
Then
\[
\frac{\partial\psi}{\partial t}=\sum_{i,j=1}^{r}\dot{X}_{i}S_{ij}W_{j}+\sum_{i,j=1}^{r}X_{i}\dot{S}_{ij}W_{j}+\sum_{i,j=1}^{r}X_{i}S_{ij}\dot{W}_{j}.
\]
It follows that
\[
\dot{S}_{ij}=\left\langle \frac{\partial\psi}{\partial t},X_{i}W_{j}\right\rangle _{\mathcal{X}\times\mathbb{S}^{2}}=\left\langle \mathcal{L}\psi,X_{i}W_{j}\right\rangle _{\mathcal{X}\times\mathbb{S}^{2}}.
\]
Therefore,
\begin{align*}
\dot{S}_{ij} & =\left\langle \mathcal{L}\sum_{i',j'=1}^{r}X_{i'}S_{i'j'}W_{j'},X_{i}W_{j}\right\rangle _{\mathcal{X}\times\mathbb{S}^{2}}\\
 & =\sum_{i',j'=1}^{r}S_{i'j'}\left\langle \mathcal{L}X_{i'}W_{j'},X_{i}W_{j}\right\rangle _{\mathcal{X}\times\mathbb{S}^{2}}.
\end{align*}

In summary,
\[
\dot{L}_{j}=\sum_{i=1}^{r}\left\langle \mathcal{L}\left(X_{i}L_{i}\right),X_{j}\right\rangle _{\mathcal{X}},
\]
\[
\dot{K}_{j}=\sum_{i=1}^{r}\left\langle \mathcal{L}\left(K_{i}W_{i}\right),W_{j}\right\rangle _{\mathbb{S}^{2}}.
\]
\[
\dot{S}_{ij}=\sum_{i',j'=1}^{r}S_{i'j'}\left\langle \mathcal{L}X_{i'}W_{j'},X_{i}W_{j}\right\rangle _{\mathcal{X}\times\mathbb{S}^{2}}.
\]

We claim that 
\[
\frac{\partial\psi}{\partial t}=\sum_{j}\dot{K}_{j}W_{j}+\sum_{i}\dot{L}_{i}X_{i}-\sum_{i,j=1}^{r}X_{i}\dot{S}_{ij}W_{j}.
\]
To show this, recall that
\[
K_{j}=\sum_{i=1}^{r}X_{i}S_{ij},\,\,\,\,\,L_{i}=\sum_{j=1}^{r}S_{ij}W_{j}.
\]
Therefore,
\[
\dot{K}_{j}=\sum_{i=1}^{r}\dot{X}_{i}S_{ij}+\sum_{i=1}^{r}X_{i}\dot{S}_{ij},
\]
and so
\[
\sum_{j}\dot{K}_{j}W_{j}=\sum_{i,j=1}^{r}\dot{X}_{i}S_{ij}W_{j}+\sum_{i,j}^{r}X_{i}\dot{S}_{ij}W_{j}.
\]
Similarly,
\[
\dot{L}_{i}=\sum_{j=1}^{r}S_{ij}\dot{W}_{j}+\sum_{i=1}^{r}\dot{S}_{ij}W_{j},
\]
and so
\[
\sum_{i}\dot{L}_{i}X_{i}=\sum_{i,j=1}^{r}X_{i}S_{ij}\dot{W}_{j}+\sum_{i,j=1}^{r}X_{i}\dot{S}_{ij}W_{j}.
\]
It follows that
\[
\sum_{j}\dot{K}_{j}W_{j}+\sum_{i}\dot{L}_{i}X_{i}=\sum_{i,j=1}^{r}\dot{X}_{i}S_{ij}W_{j}+2\sum_{i,j=1}^{r}X_{i}\dot{S}_{ij}W_{j}+\sum_{i,j=1}^{r}X_{i}S_{ij}\dot{W}_{j}.
\]
Therefore,
\begin{align*}
\sum_{j}\dot{K}_{j}W_{j}+\sum_{i}\dot{L}_{i}X_{i}-\sum_{i,j=1}^{r}X_{i}\dot{S}_{ij}W_{j} & =\sum_{i,j=1}^{r}\dot{X}_{i}S_{ij}W_{j}+\sum_{i,j=1}^{r}X_{i}\dot{S}_{ij}W_{j}+\sum_{i,j=1}^{r}X_{i}S_{ij}\dot{W}_{j}\\
 & =\frac{\partial}{\partial t}\left(\sum_{i,j=1}^{r}X_{i}S_{ij}W_{j}\right).
\end{align*}

\subsection{Derivation of DSA for collocation-based DLR \label{subsec:Derivation-of-DSA}}

\subsubsection{Diffusion limit}

Write equation (\ref{eq:L step}) as 
\[
\left(I+\varepsilon\sum_{k=1}^{3}\Omega_{k}M_{t}^{-1}A_{k}\right)\mathbf{L}=\frac{1}{4\pi}\left(I-\varepsilon^{2}M_{t}^{-1}M_{a}\right)\int_{\mathbb{S}^{2}}\mathbf{L}d\boldsymbol{\Omega}+\varepsilon M_{t}^{-1}\mathbf{Q}.
\]
In the diffusion limit,
\begin{align*}
\left(I+\varepsilon\sum_{k=1}^{3}\Omega_{k}M_{t}^{-1}A_{k}\right)^{-1} & =I-\varepsilon\sum_{k=1}^{3}\Omega_{k}M_{t}^{-1}A_{k}+\\
 & \,\,\,\,\,\,\varepsilon^{2}\left(\sum_{k=1}^{3}\Omega_{k}M_{t}^{-1}A_{k}\right)^{2}+\mathcal{O}\left(\varepsilon^{3}\right).
\end{align*}
To order $\mathcal{O}\left(\varepsilon^{2}\right)$, $L_{j}\left(\boldsymbol{\Omega}\right)$
can be written as linear combinations of $1$, $\Omega_{k}$, $\Omega_{k}\Omega_{l}$,
\[
L_{j}\left(\boldsymbol{\Omega}\right)=a_{j}+\sum_{k=1}^{3}a_{j,k}\Omega_{k}+\sum_{k,l=1}^{3}a_{j,k,l}\Omega_{k}\Omega_{l}+\mathcal{O}\left(\varepsilon^{3}\right).
\]
It follows that there are coefficients $b_{d'}$, $b_{d',k}$, $b_{d',k,l}$
for which
\begin{equation}
W_{d'}\left(\boldsymbol{\Omega}\right)=b_{d'}+\sum_{k=1}^{3}b_{d',k}\left(\boldsymbol{\Omega}\right)_{k}+\sum_{k,l=1}^{3}b_{d',k,l}\left(\boldsymbol{\Omega}\right)_{k}\left(\boldsymbol{\Omega}\right)_{l}+\mathcal{O}\left(\varepsilon^{3}\right),\,\,\,\,\,d'=1,\ldots,r_{0}.\label{eq:W_d' in terms of polynomials}
\end{equation}

Assume that we can solve, up to errors of order $\mathcal{O}\left(\varepsilon^{3}\right)$,
equations (\ref{eq:W_d' in terms of polynomials}) for the basis functions
$1$, $\left(\boldsymbol{\Omega}\right)_{k}\left(\boldsymbol{\Omega}\right)_{l}$,
$\left(\boldsymbol{\Omega}\right)_{k}\left(\boldsymbol{\Omega}\right)_{l}$,
$k,l=1,2,3$, in terms of $W_{d'}\left(\boldsymbol{\Omega}\right)$,
$d'=1,\ldots,r_{0}$. Then there are coefficients $\alpha_{d'}^{\left(i,j\right)}$
for which
\[
\left(\boldsymbol{\Omega}\right)_{i}\left(\boldsymbol{\Omega}\right)_{j}=\sum_{d'=1}^{r_{0}}\alpha_{d'}^{\left(i,j\right)}W_{d'}\left(\boldsymbol{\Omega}\right)+\mathcal{O}\left(\varepsilon^{3}\right).
\]
Therefore,
\[
\left(\boldsymbol{\Omega}_{d}\right)_{i}\left(\boldsymbol{\Omega}_{d}\right)_{j}=\sum_{d'=1}^{r'}\alpha_{d'}^{\left(i,j\right)}W_{d'}\left(\boldsymbol{\Omega}_{d}\right)+\mathcal{O}\left(\varepsilon^{3}\right),\,\,\,\,\,d=1,\ldots,r_{0}.
\]
It follows that, for each $i,j=1,2,3$,
\[
\boldsymbol{\alpha}^{\left(i,j\right)}=\hat{W}^{-1}\left[\left(\boldsymbol{\Omega}_{d}\right)_{i}\left(\boldsymbol{\Omega}_{d}\right)_{j}\right]_{d=1}^{r_{0}}+\mathcal{O}\left(\varepsilon^{3}\right).
\]
We then calculate that
\begin{align*}
\int_{\mathbb{S}^{2}}\left(\boldsymbol{\Omega}\right)_{i}\left(\boldsymbol{\Omega}\right)_{j}d\boldsymbol{\Omega} & =\sum_{d'}\alpha_{d'}^{\left(i,j\right)}\int_{\mathbb{S}^{2}}W_{d'}\left(\boldsymbol{\Omega}\right)d\boldsymbol{\Omega}+\mathcal{O}\left(\varepsilon^{3}\right)\\
 & =\boldsymbol{\beta}^{\text{T}}\boldsymbol{\alpha}^{\left(i,j\right)}+\mathcal{O}\left(\varepsilon^{3}\right)\\
 & =\boldsymbol{\beta}^{\text{T}}\hat{W}^{-1}\left[\left(\boldsymbol{\Omega}_{d}\right)_{i}\left(\boldsymbol{\Omega}_{d}\right)_{j}\right]_{d=1}^{r}+\mathcal{O}\left(\varepsilon^{3}\right).
\end{align*}

\subsubsection{Diffusion Synthetic Acceleration (DSA) for SN-like DLR method}

Recall that
\[
\left(\boldsymbol{\Omega}_{d}\cdot\mathbf{G}+F_{d}\right)\boldsymbol{\psi}_{d}+M_{t}\boldsymbol{\psi}_{d}=\frac{M_{s}}{4\pi}\boldsymbol{\phi}+MQ_{d},
\]
\[
\boldsymbol{\phi}=\boldsymbol{\beta}^{\text{T}}\left(\hat{W}^{-1}\boldsymbol{\psi}\right).
\]
Here $\hat{W}=\left[W_{i}\left(\boldsymbol{\Omega}_{j}\right)\right]_{i,j=1}^{r}$,
$\beta_{i}=\int_{\mathbb{S}^{2}}W_{i}\left(\boldsymbol{\Omega}\right)d\boldsymbol{\Omega}$,
and 
\begin{align*}
\mathbf{v}^{\text{T}}F_{d}\mathbf{u} & =-\sum_{\Gamma\in\mathcal{F}}\int_{\Gamma}\boldsymbol{\Omega}_{d}\cdot\mathbf{n}\left\llbracket u\right\rrbracket \left\{ v\right\} dS+\frac{1}{2}\sum_{\Gamma\in\mathcal{F}}\int_{\Gamma}\left|\boldsymbol{\Omega}_{d}\cdot\mathbf{n}\right|\left\llbracket u\right\rrbracket \left\llbracket v\right\rrbracket dS\\
 & =\mathbf{v}^{\text{T}}\left(\boldsymbol{\Omega}_{d}\cdot\mathbf{F}\right)\mathbf{u}+\mathbf{v}^{\text{T}}F_{d}^{\left\llbracket \right\rrbracket }\mathbf{u}.
\end{align*}
Let $\hat{\boldsymbol{\psi}}_{d}$ denote the solution to 
\[
\left(-\boldsymbol{\Omega}_{d}\cdot\mathbf{G}+F_{d}\right)\hat{\boldsymbol{\psi}}_{d}+M_{t}\boldsymbol{\psi}_{d}=\frac{M_{s}}{4\pi}\boldsymbol{\phi}+MQ_{d},
\]
\[
\boldsymbol{\phi}=\boldsymbol{\beta}^{\text{T}}\left(\hat{W}^{-1}\boldsymbol{\psi}\right).
\]
Here 
\begin{align*}
\mathbf{v}^{\text{T}}F_{d}\mathbf{u} & =-\sum_{\Gamma\in\mathcal{F}}\int_{\Gamma}\boldsymbol{\Omega}_{d}\cdot\mathbf{n}\left\llbracket u\right\rrbracket \left\{ v\right\} dS+\frac{1}{2}\sum_{\Gamma\in\mathcal{F}}\int_{\Gamma}\left|\boldsymbol{\Omega}_{d}\cdot\mathbf{n}\right|\left\llbracket u\right\rrbracket \left\llbracket v\right\rrbracket dS\\
 & =\mathbf{v}^{\text{T}}\left(\boldsymbol{\Omega}_{d}\cdot\mathbf{F}\right)\mathbf{u}+\mathbf{v}^{\text{T}}F_{d}^{\left\llbracket \right\rrbracket }\mathbf{u}.
\end{align*}

Now, 
\begin{align*}
\mathbf{v}^{\text{T}}F_{\mathcal{R}\left(d\right)}\mathbf{u} & =-\sum_{\Gamma\in\mathcal{F}}\int_{\Gamma}\boldsymbol{\Omega}_{\mathcal{R}\left(d\right)}\cdot\mathbf{n}\left\llbracket u\right\rrbracket \left\{ v\right\} dS+\frac{1}{2}\sum_{\Gamma\in\mathcal{F}}\int_{\Gamma}\left|\boldsymbol{\Omega}_{\mathcal{R}\left(d\right)}\cdot\mathbf{n}\right|\left\llbracket u\right\rrbracket \left\llbracket v\right\rrbracket dS\\
 & =-\mathbf{v}^{\text{T}}\left(\boldsymbol{\Omega}_{d}\cdot\mathbf{F}\right)\mathbf{u}+\mathbf{v}^{\text{T}}F_{d}^{\left\llbracket \right\rrbracket }\mathbf{u}.
\end{align*}
It follows that
\[
F_{\mathcal{R}\left(d\right)}\boldsymbol{\psi}_{\mathcal{R}\left(d\right)}=-\left(\boldsymbol{\Omega}_{d}\cdot\mathbf{F}\right)\boldsymbol{\psi}_{\mathcal{R}\left(d\right)}+F_{d}^{\left\llbracket \right\rrbracket }\boldsymbol{\psi}_{\mathcal{R}\left(d\right)},
\]
\[
F_{d}\boldsymbol{\psi}_{d}=\left(\boldsymbol{\Omega}_{d}\cdot\mathbf{F}\right)\boldsymbol{\psi}_{d}+F_{d}^{\left\llbracket \right\rrbracket }\boldsymbol{\psi}_{d}.
\]

Now, define
\[
\boldsymbol{\psi}_{d}^{+}=\frac{\boldsymbol{\psi}_{\mathcal{R}\left(d\right)}+\boldsymbol{\psi}_{d}}{2},\,\,\,\,\,\boldsymbol{\psi}_{d}^{-}=\frac{\boldsymbol{\psi}_{d}-\boldsymbol{\psi}_{\mathcal{R}\left(d\right)}}{2}.
\]
We compute that
\[
\frac{1}{2}\left(F_{\mathcal{R}\left(d\right)}\boldsymbol{\psi}_{\mathcal{R}\left(d\right)}+F_{d}\boldsymbol{\psi}_{d}\right)=\left(\boldsymbol{\Omega}_{d}\cdot\mathbf{F}\right)\boldsymbol{\psi}_{d}^{-}+F_{d}^{\left\llbracket \right\rrbracket }\boldsymbol{\psi}_{d}^{+},
\]
\[
\frac{1}{2}\left(F_{d}\boldsymbol{\psi}_{d}-F_{\mathcal{R}\left(d\right)}\boldsymbol{\psi}_{\mathcal{R}\left(d\right)}\right)=\left(\boldsymbol{\Omega}_{d}\cdot\mathbf{F}\right)\boldsymbol{\psi}_{d}^{+}+F_{d}^{\left\llbracket \right\rrbracket }\boldsymbol{\psi}_{d}^{-}.
\]
Now,
\[
\left(-\boldsymbol{\Omega}_{d}\cdot\mathbf{G}+F_{\mathcal{R}\left(d\right)}\right)\boldsymbol{\psi}_{\mathcal{R}\left(d\right)}+M_{t}\boldsymbol{\psi}_{\mathcal{R}\left(d\right)}=\frac{M_{s}}{4\pi}\boldsymbol{\phi}+MQ_{\mathcal{R}\left(d\right)}.
\]
\[
\left(\boldsymbol{\Omega}_{d}\cdot\mathbf{G}+F_{d}\right)\boldsymbol{\psi}_{d}+M_{t}\boldsymbol{\psi}_{d}=\frac{M_{s}}{4\pi}\boldsymbol{\phi}+MQ_{d}.
\]
Add and subtract these equations and divide by $2$:
\[
\boldsymbol{\Omega}_{d}\cdot\left(\mathbf{G}+\mathbf{F}\right)\boldsymbol{\psi}_{d}^{-}+F_{d}^{\left\llbracket \right\rrbracket }\boldsymbol{\psi}_{d}^{+}+M_{t}\boldsymbol{\psi}_{d}^{+}=\frac{M_{s}}{4\pi}\boldsymbol{\phi}+MQ_{d}^{+},
\]
\[
\boldsymbol{\Omega}_{d}\cdot\left(\mathbf{G}+\mathbf{F}\right)\boldsymbol{\psi}_{d}^{+}+F_{d}^{\left\llbracket \right\rrbracket }\boldsymbol{\psi}_{d}^{-}+M_{t}\boldsymbol{\psi}_{d}^{-}=MQ_{d}^{-}.
\]

Now, define 
\[
\mathbf{v}^{\text{T}}\tilde{F}_{d}\mathbf{u}=\sum_{\Gamma\in\mathcal{F}}\int_{\Gamma}\boldsymbol{\Omega}_{d}\cdot\mathbf{n}\left\llbracket v\right\rrbracket \left\{ u\right\} dS+\frac{1}{2}\sum_{\Gamma\in\mathcal{F}}\int_{\Gamma}\left|\boldsymbol{\Omega}_{d}\cdot\mathbf{n}\right|\left\llbracket u\right\rrbracket \left\llbracket v\right\rrbracket dS.
\]
It turns out that
\[
\boldsymbol{\Omega}_{d}\cdot\mathbf{G}+F_{d}=-\boldsymbol{\Omega}_{d}\cdot\mathbf{G}^{\text{T}}+\tilde{F}_{d}.
\]
Recall that 
\begin{align*}
\mathbf{v}^{\text{T}}F_{d}\mathbf{u} & =-\sum_{\Gamma\in\mathcal{F}}\int_{\Gamma}\boldsymbol{\Omega}_{d}\cdot\mathbf{n}\left\llbracket u\right\rrbracket \left\{ v\right\} dS+\frac{1}{2}\sum_{\Gamma\in\mathcal{F}}\int_{\Gamma}\left|\boldsymbol{\Omega}_{d}\cdot\mathbf{n}\right|\left\llbracket u\right\rrbracket \left\llbracket v\right\rrbracket dS\\
 & =\mathbf{v}^{\text{T}}\left(\boldsymbol{\Omega}_{d}\cdot\mathbf{F}\right)\mathbf{u}+\mathbf{v}^{\text{T}}F_{d}^{\left\llbracket \right\rrbracket }\mathbf{u}.
\end{align*}
Also,
\[
\tilde{F}_{d}=-\boldsymbol{\Omega}_{d}\cdot\mathbf{F}^{\text{T}}+F_{d}^{\left\llbracket \right\rrbracket }.
\]
It follows that 
\[
\boldsymbol{\Omega}_{d}\cdot\mathbf{G}+\boldsymbol{\Omega}_{d}\cdot\mathbf{F}+F_{d}^{\left\llbracket \right\rrbracket }=-\boldsymbol{\Omega}_{d}\cdot\mathbf{G}^{\text{T}}-\boldsymbol{\Omega}_{d}\cdot\mathbf{F}^{\text{T}}+F_{d}^{\left\llbracket \right\rrbracket }.
\]
Therefore, we obtain the equations for the even and odd components,
\[
-\boldsymbol{\Omega}_{d}\cdot\left(\mathbf{G}+\mathbf{F}\right)^{\text{T}}\boldsymbol{\psi}_{d}^{-}+F_{d}^{\left\llbracket \right\rrbracket }\boldsymbol{\psi}_{d}^{+}+M_{t}\boldsymbol{\psi}_{d}^{+}=\frac{M_{s}}{4\pi}\boldsymbol{\phi}+MQ_{d}^{+},
\]
\[
\boldsymbol{\Omega}_{d}\cdot\left(\mathbf{G}+\mathbf{F}\right)\boldsymbol{\psi}_{d}^{+}+F_{d}^{\left\llbracket \right\rrbracket }\boldsymbol{\psi}_{d}^{-}+M_{t}\boldsymbol{\psi}_{d}^{-}=MQ_{d}^{-}.
\]

Now,
\begin{align*}
\boldsymbol{\psi}_{d}^{-} & =\left(I+M_{t}^{-1}F_{d}^{\left\llbracket \right\rrbracket }\right)^{-1}M_{t}^{-1}\left(MQ_{d}^{-}\right)-\left(I+M_{t}^{-1}F_{d}^{\left\llbracket \right\rrbracket }\right)^{-1}M_{t}^{-1}\boldsymbol{\Omega}_{d}\cdot\left(\mathbf{G}+\mathbf{F}\right)\boldsymbol{\psi}_{d}^{+}\\
 & =-M_{t}^{-1}\boldsymbol{\Omega}_{d}\cdot\left(\mathbf{G}+\mathbf{F}\right)\boldsymbol{\psi}_{d}^{+}+M_{t}^{-1}\left(MQ_{d}^{-}\right)+\\
 & \,\,\,\,\,\,\,\,\,\,\,\,\,\left(\left(I+M_{t}^{-1}F_{d}^{\left\llbracket \right\rrbracket }\right)^{-1}-I\right)M_{t}^{-1}\left(MQ_{d}^{-}\right)-\\
 & \,\,\,\,\,\,\,\,\,\,\,\,\,\,\left(\left(I+M_{t}^{-1}F_{d}^{\left\llbracket \right\rrbracket }\right)^{-1}-I\right)M_{t}^{-1}\boldsymbol{\Omega}_{d}\cdot\left(\mathbf{G}+\mathbf{F}\right)\boldsymbol{\psi}_{d}^{+}.
\end{align*}
Therefore, to leading order
\begin{align*}
-\boldsymbol{\Omega}_{d}\cdot\left(\mathbf{G}+\mathbf{F}\right)^{\text{T}}\left(-M_{t}^{-1}\boldsymbol{\Omega}_{d}\cdot\left(\mathbf{G}+\mathbf{F}\right)\boldsymbol{\psi}_{d}^{+}\right)+F_{d}^{\left\llbracket \right\rrbracket }\boldsymbol{\psi}_{d}^{+}+M_{t}\boldsymbol{\psi}_{d}^{+} & \approx\frac{M_{s}}{4\pi}\boldsymbol{\phi}+MQ_{d}^{+}+\\
 & \,\,\,\,\boldsymbol{\Omega}_{d}\cdot\left(\mathbf{G}+\mathbf{F}\right)^{\text{T}}M_{t}^{-1}\left(MQ_{d}^{-}\right).
\end{align*}
Also, to leading order,
\[
\boldsymbol{\psi}_{d}=M_{t}^{-1}\frac{M_{s}}{4\pi}\boldsymbol{\phi}+\mathcal{O}\left(\varepsilon\right)=\frac{1}{4\pi}\boldsymbol{\phi}+\mathcal{O}\left(\varepsilon\right),
\]
and so
\[
\boldsymbol{\psi}_{d}^{+}=\frac{1}{4\pi}\boldsymbol{\phi}+\mathcal{O}\left(\varepsilon\right),\,\,\,\,\,\,\boldsymbol{\psi}_{d}^{-}=\mathcal{O}\left(\varepsilon\right).
\]
We can drop the term $\boldsymbol{\Omega}_{d}\cdot\left(\mathbf{G}+\mathbf{F}\right)^{\text{T}}M_{t}^{-1}\left(MQ_{d}^{-}\right)$
since we assume that $Q_{d}=\mathcal{O}\left(\varepsilon\right)$.
It follows that
\begin{equation}
-\boldsymbol{\Omega}_{d}\cdot\left(\mathbf{G}+\mathbf{F}\right)^{\text{T}}\left(-M_{t}^{-1}\boldsymbol{\Omega}_{d}\cdot\left(\mathbf{G}+\mathbf{F}\right)\frac{1}{4\pi}\boldsymbol{\phi}\right)+F_{d}^{\left\llbracket \right\rrbracket }\frac{1}{4\pi}\boldsymbol{\phi}+M_{t}\frac{1}{4\pi}\boldsymbol{\phi}\approx\frac{M_{s}}{4\pi}\boldsymbol{\phi}+MQ_{d}^{+}.\label{eq:phi moment system}
\end{equation}

Now consider 
\[
-\boldsymbol{\Omega}_{d}\cdot\left(\mathbf{G}+\mathbf{F}\right)^{\text{T}}\left(-M_{t}^{-1}\boldsymbol{\Omega}_{d}\cdot\left(\mathbf{G}+\mathbf{F}\right)\frac{1}{4\pi}\boldsymbol{\phi}\right)+F_{d}^{\left\llbracket \right\rrbracket }\frac{1}{4\pi}\boldsymbol{\phi}+M_{a}\frac{1}{4\pi}\boldsymbol{\phi}\approx MQ_{d}^{+}.
\]
Define 
\[
\boldsymbol{\alpha}=\hat{W}^{-1}\left[1\right]_{d=1}^{r},
\]
\[
\boldsymbol{\alpha}^{\left(i,j\right)}=\hat{W}^{-1}\left[\left(\boldsymbol{\Omega}_{d}\right)_{i}\left(\boldsymbol{\Omega}_{d}\right)_{j}\right]_{d=1}^{r},
\]
\[
\gamma=\boldsymbol{\beta}^{\text{T}}\boldsymbol{\alpha},
\]
\[
\gamma^{\left(i,j\right)}=\boldsymbol{\beta}^{\text{T}}\boldsymbol{\alpha}^{\left(i,j\right)}.
\]
Let 
\[
\tilde{\boldsymbol{\phi}}_{d}=\left(\boldsymbol{\Omega}_{d}\right)_{i}\left(\boldsymbol{\Omega}_{d}\right)_{j}\left(\mathbf{G}_{i}+\mathbf{F}_{i}\right)^{\text{T}}M_{t}^{-1}\left(\mathbf{G}_{j}+\mathbf{F}_{j}\right)\boldsymbol{\phi}.
\]
Then 
\begin{align*}
\boldsymbol{\beta}^{\text{T}}\hat{W}^{-1}\left[\tilde{\boldsymbol{\phi}}_{d}\right]_{d} & =\boldsymbol{\beta}^{\text{T}}\hat{W}^{-1}\left[\left(\boldsymbol{\Omega}_{d}\right)_{i}\left(\boldsymbol{\Omega}_{d}\right)_{j}\right]_{d=1}^{r}\left(\mathbf{G}_{i}+\mathbf{F}_{i}\right)^{\text{T}}M_{t}^{-1}\left(\mathbf{G}_{j}+\mathbf{F}_{j}\right)\boldsymbol{\phi}\\
 & =\boldsymbol{\beta}^{\text{T}}\boldsymbol{\alpha}^{\left(i,j\right)}\left(\mathbf{G}_{i}+\mathbf{F}_{i}\right)^{\text{T}}M_{t}^{-1}\left(\mathbf{G}_{j}+\mathbf{F}_{j}\right)\boldsymbol{\phi}.
\end{align*}
It follows that the DSA matrix can be written as
\[
\sum_{i,j=1}^{3}\gamma^{\left(i,j\right)}\left(\mathbf{G}_{i}+\mathbf{F}_{i}\right)^{\text{T}}M_{t}^{-1}\left(\mathbf{G}_{j}+\mathbf{F}_{j}\right)\boldsymbol{\phi}+F_{0}\frac{1}{4\pi}\boldsymbol{\phi}+\gamma M_{a}\frac{1}{4\pi}\boldsymbol{\phi}.
\]
Here
\[
\mathbf{v}^{\text{T}}F_{0}\mathbf{u}=\frac{1}{4\pi}\sum_{\Gamma\in\mathcal{F}}\left(\frac{1}{2}\int_{\Gamma}\boldsymbol{\beta}^{\text{T}}\hat{W}^{-1}\left[\left|\boldsymbol{\Omega}_{d}\cdot\mathbf{n}\right|\right]_{d=1}^{r}\right)\left\llbracket u\right\rrbracket \left\llbracket v\right\rrbracket dS.
\]

Now, suppose that
\[
\left(\boldsymbol{\Omega}_{d}\right)_{i}\left(\boldsymbol{\Omega}_{d}\right)_{j}=\sum_{d'=1}^{r'}\alpha_{d'}^{\left(i,j\right)}W_{d'}\left(\boldsymbol{\Omega}_{d}\right),
\]
so that 
\[
\boldsymbol{\alpha}^{\left(i,j\right)}=\hat{W}^{-1}\left[\left(\boldsymbol{\Omega}_{d}\right)_{i}\left(\boldsymbol{\Omega}_{d}\right)_{j}\right]_{d=1}^{r}.
\]
Also,
\[
\int_{\mathbb{S}^{2}}\left(\boldsymbol{\Omega}\right)_{i}\left(\boldsymbol{\Omega}\right)_{j}d\boldsymbol{\Omega}=\sum_{d'}\alpha_{d'}^{\left(i,j\right)}\int_{\mathbb{S}^{2}}W_{d'}\left(\boldsymbol{\Omega}\right)d\boldsymbol{\Omega}=\boldsymbol{\beta}^{\text{T}}\boldsymbol{\alpha}^{\left(i,j\right)}=\boldsymbol{\beta}^{\text{T}}\hat{W}^{-1}\left[\left(\boldsymbol{\Omega}_{d}\right)_{i}\left(\boldsymbol{\Omega}_{d}\right)_{j}\right]_{d=1}^{r}.
\]

Now,
\[
\boldsymbol{\Omega}_{d}\cdot\left(\mathbf{G}+\mathbf{F}\right)^{\text{T}}M_{t}^{-1}\boldsymbol{\Omega}_{d}\cdot\left(\mathbf{G}+\mathbf{F}\right)=\sum_{i,j=1}^{3}\left(\boldsymbol{\Omega}_{d}\right)_{i}\left(\boldsymbol{\Omega}_{d}\right)_{j}\left(\mathbf{G}_{i}+\mathbf{F}_{i}\right)^{\text{T}}M_{t}^{-1}\left(\mathbf{G}_{j}+\mathbf{F}_{j}\right).
\]
Let 
\[
\tilde{\boldsymbol{\phi}}_{d}=\left(\boldsymbol{\Omega}_{d}\right)_{i}\left(\boldsymbol{\Omega}_{d}\right)_{j}\left(\mathbf{G}_{i}+\mathbf{F}_{i}\right)^{\text{T}}M_{t}^{-1}\left(\mathbf{G}_{j}+\mathbf{F}_{j}\right)\boldsymbol{\phi}.
\]
Then 
\begin{align*}
\boldsymbol{\beta}^{\text{T}}\hat{W}^{-1}\left[\tilde{\boldsymbol{\phi}}_{d}\right]_{d} & =\boldsymbol{\beta}^{\text{T}}\hat{W}^{-1}\left[\left(\boldsymbol{\Omega}_{d}\right)_{i}\left(\boldsymbol{\Omega}_{d}\right)_{j}\right]_{d=1}^{r}\left(\mathbf{G}_{i}+\mathbf{F}_{i}\right)^{\text{T}}M_{t}^{-1}\left(\mathbf{G}_{j}+\mathbf{F}_{j}\right)\boldsymbol{\phi}\\
 & =\left(\int_{\mathbb{S}^{2}}\left(\boldsymbol{\Omega}\right)_{i}\left(\boldsymbol{\Omega}\right)_{j}d\boldsymbol{\Omega}\right)\left(\mathbf{G}_{i}+\mathbf{F}_{i}\right)^{\text{T}}M_{t}^{-1}\left(\mathbf{G}_{j}+\mathbf{F}_{j}\right)\boldsymbol{\phi}\\
 & =\frac{4\pi}{3}\delta_{i.j}\left(\mathbf{G}_{i}+\mathbf{F}_{i}\right)^{\text{T}}M_{t}^{-1}\left(\mathbf{G}_{j}+\mathbf{F}_{j}\right)\boldsymbol{\phi}.
\end{align*}
Applying $\boldsymbol{\beta}^{\text{T}}\hat{W}^{-1}$ to equation
(\ref{eq:phi moment system}), it follows that

\[
\frac{1}{3}\left(\mathbf{G}+\mathbf{F}\right)^{\text{T}}M_{t}^{-1}\cdot\left(\mathbf{G}+\mathbf{F}\right)\boldsymbol{\phi}+F_{0}\boldsymbol{\phi}+M_{a}\boldsymbol{\phi}\approx MQ^{+}.
\]
Here
\[
\mathbf{v}^{\text{T}}F_{0}\mathbf{u}=\frac{1}{4\pi}\sum_{\Gamma\in\mathcal{F}}\left(\frac{1}{2}\int_{\Gamma}\boldsymbol{\beta}^{\text{T}}\hat{W}^{-1}\left[\left|\boldsymbol{\Omega}_{d}\cdot\mathbf{n}\right|\right]_{d=1}^{r}\right)\left\llbracket u\right\rrbracket \left\llbracket v\right\rrbracket dS.
\]
We also used that
\[
\boldsymbol{\beta}^{\text{T}}\hat{W}^{-1}\left[1\right]_{d=1}^{r}=4\pi+\mathcal{O}\left(\varepsilon\right),
\]
\[
\boldsymbol{\beta}^{\text{T}}\hat{W}^{-1}\left[\left(\boldsymbol{\Omega}_{d}\right)_{i}\left(\boldsymbol{\Omega}_{d}\right)_{j}\right]_{d=1}^{r}=\int_{\mathbb{S}^{2}}\left(\boldsymbol{\Omega}\right)_{i}\left(\boldsymbol{\Omega}\right)_{j}d\boldsymbol{\Omega}+\mathcal{O}\left(\varepsilon\right).
\]

~~~~~~~

\section{Acknowledgements}

This work was performed under the auspices of the U.S. Department
of Energy by Lawrence Livermore National Laboratory under Contract
DE-AC52-07NA27344.
This material is based upon work supported by the U.S. Department of Energy, Office of Science, Office of Mathematical Multifaceted Integrated Capability Center (MMICC) under Award Number DE-SC-000XXXX.

This work has been reviewed for unlimited public release as LLNL-JRNL-2015102.

\bibliographystyle{plain}
\addcontentsline{toc}{section}{\refname}\bibliography{refs}

\end{document}